\documentclass[11pt]{amsart}

\usepackage[T1]{fontenc}
\usepackage{lmodern}
\usepackage{amsmath,amssymb}
\usepackage{microtype}
\usepackage{needspace}
\usepackage[margin=1.15in]{geometry}

\numberwithin{equation}{section}
\newtheorem{theorem}{Theorem}[section]
\newtheorem{proposition}[theorem]{Proposition}
\newtheorem{corollary}[theorem]{Corollary}
\newtheorem{lemma}[theorem]{Lemma}
\theoremstyle{definition}
\newtheorem{example}[theorem]{Example}
\theoremstyle{remark}
\newtheorem{remark}[theorem]{Remark}

\title[Zero-cycles on quotients of multiple double covers]
{Zero-cycles on quotients of multiple double covers and their products}
\author{Wenchuan Hu}
\address{School of Mathematics, Sichuan University, Chengdu 610064, P. R. China}
\email{huwenchuan@gmail.com}
\email{wenchuan@scu.edu.cn}
\date{\today}

\begin{document}
\raggedbottom

\begin{abstract}
We prove the 0-cycle form of the generalized Bloch conjecture for
smooth projective models of subgroup quotients of multiple double covers of
projective space, with smooth branch hypersurfaces and reduced simple
normal crossings branch divisor. Our integral character criterion yields,
in every dimension at least four, simply connected varieties of general
type with 0-cycle Chow group \(\mathbb{Z}\), canonical volume four,
and no birational product decomposition. We also determine optimal
integral supports for componentwise product quotients and character
decompositions of degree zero Chow groups. Within this covering class,
we classify the minimum volume models and compute their pluricanonical map degrees.
\end{abstract}

\maketitle

\section{Introduction}

The generalized Bloch conjecture relates the vanishing of holomorphic
forms to the support of 0-cycles. We establish its 0-cycle form
for finite quotients of multiple double covers and use the resulting
integral Chow criterion to construct nonproduct varieties of general type
in every dimension at least four.

In this paper, varieties are reduced and irreducible over \(\mathbb{C}\),
and Chow groups have integral coefficients. Supports may be reducible.
For a projective variety or a closed reduced projective subscheme \(V\),
write \(\operatorname{CH}_0(V)\) for its group of 0-cycles modulo
rational equivalence, and define
\begin{equation*}
A_0(V)=\ker\bigl(\deg:\operatorname{CH}_0(V)\longrightarrow\mathbb{Z}\bigr).
\end{equation*}
The class of a closed point is denoted by \([x]\). For \(V\) smooth,
let \(H^{p,0}(V)=H^0(V,\Omega_V^p)\) and
\(h^{p,0}(V)=\dim_{\mathbb{C}}H^{p,0}(V)\), where \(\Omega_V^p\)
is the sheaf of regular differential \(p\)-forms. If \(\dim V=m\),
write \(\omega_V=\Omega_V^m\), \(K_V\) for a canonical divisor, and
\(p_g(V)=h^{m,0}(V)\). The symbols \(\sim\) and
\(\sim_{\mathbb{Q}}\) denote linear and rational linear equivalence.

To formulate the integral support problem, let \(V\) be a smooth
projective variety and define
\begin{equation*}
s_0(V)=\min\left\{\dim Z:
\operatorname{CH}_0(Z)\longrightarrow\operatorname{CH}_0(V)
\text{ is surjective}\right\},
\end{equation*}
where \(Z\) ranges over nonempty closed reduced subsets and the map is
induced by inclusion. The 0-cycle form of the generalized Bloch
conjecture asserts, for every integer \(s\geq0\), that
\begin{equation*}
H^{p,0}(V)=0\quad(p>s)\qquad\Longrightarrow\qquad s_0(V)\leq s;
\end{equation*}
see \cite[Conjecture~1]{Hu} and \cite[Conjecture~1.11]{Voi}.
The reverse implication is the generalized Mumford theorem
\cite{Mum,BS}, \cite[Theorem~3.13]{Voi}.

Let \(X\longrightarrow\mathbb{P}^n\), \(n\geq2\), be the fiber
product of \(r\) double covers with smooth branch hypersurfaces of
degrees \(2d_i\), where \(d_i\geq1\), whose sum is reduced with
simple normal crossings. Let \(\tau_i\) be the involution changing the sign of the
\(i\)-th root and \(G=\langle\tau_1,\ldots,\tau_r\rangle\) the
covering group. For \(I\subseteq\{1,\ldots,r\}\), define
\(d_I=\sum_{i\in I}d_i\), and let \(\chi_I\) be the character
with \(\chi_I(\tau_i)=-1\) for \(i\in I\) and
\(\chi_I(\tau_i)=1\) otherwise. For every subgroup \(H\subseteq G\)
and every smooth projective model \(Y\) of \(X/H\),
Theorem~\ref{theorem:single-cover} proves
\[
\deg:\operatorname{CH}_0(Y)\xrightarrow{\ \cong\ }\mathbb{Z}
\quad\Longleftrightarrow\quad
d_I\leq n\text{ whenever }\chi_I|_H=1.
\]
The same condition is equivalent to the vanishing of all positive-degree
holomorphic forms on \(Y\). Only the characters trivial on \(H\)
enter the criterion; no bound on the total branch degree is assumed.
For componentwise subgroup quotients of products,
Theorem~\ref{theorem:products} establishes the full support formula
\[
s_0(Y)=\max\{p\geq0:H^{p,0}(Y)\neq0\}.
\]
Diagonal subgroups and positive dimensional minimal supports are included.

The integral criterion yields a uniform family of nonproduct varieties
of general type, as described in Theorem~\ref{proposition:general type}. For \(n\geq4\), choose three general branch
hypersurfaces of degrees \(4,4,2n-4\), so that
\((d_1,d_2,d_3)=(2,2,n-2)\), and let
\(H=\langle\tau_1\tau_2\tau_3\rangle\). The nonempty subsets
whose characters are trivial on \(H\) are exactly the three pairs,
with degree sums
\[
d_{\{1,2\}}=4,\qquad d_{\{1,3\}}=d_{\{2,3\}}=n.
\]
Theorem~\ref{theorem:single-cover} therefore gives
\(\operatorname{CH}_0(Y)\cong\mathbb{Z}\) by the degree map.
On the other hand, the total half-branch degree is \(n+2\), and
\(H\) contains no individual coordinate involution. Thus the quotient
\(Q=X/H\) has canonical singularities and
\(K_Q\sim_{\mathbb{Q}}\overline f^{\,*}A\), where
\(\overline f:Q\longrightarrow\mathbb{P}^n\) is the induced finite
morphism and \(A\) is a hyperplane divisor. Its canonical divisor is
ample, so \(Y\) is of general type. Moreover, the canonical model
\(Q\) has Picard number one, excluding any birational decomposition
of \(Y\) into two positive dimensional smooth projective factors.
The construction is a fiber product over a common base, not a Cartesian
product. The general type conclusion and integral point support are
therefore obtained from the same branch degrees and subgroup characters.

For the varieties furnished by Theorem~\ref{proposition:general type},
the fundamental group, canonical index, and pluricanonical results give
further geometric information; they impose no additional hypotheses
on the Chow criterion. The fixed point calculation proves
\(\pi_1(Y)=1\), and the canonical divisor calculation gives
\(\operatorname{vol}(K_Y)=4\). Here \(\pi_1\) is the topological
fundamental group and \(\operatorname{vol}(K_Y)\) the canonical volume.
More generally, among general type models in this covering class with
\(A_0=0\), the smooth models are simply connected and their canonical models have
index two: \(K_Q\) is not Cartier, but \(2K_Q\) is
Cartier. Theorem~\ref{theorem:minimum volume} proves that four is the
minimum canonical volume and classifies equality within this class.
The pluricanonical map formulas then distinguish the birational geometry
of models with the same integral Chow group and canonical volume.
Neither the volume bound nor the dimension-four threshold is asserted
for arbitrary varieties with \(A_0=0\).

Relevant precedents include the odd-dimensional Kummer resolutions of
Bloch--Srinivas \cite{BS}, with divisor support, and the Calabi--Yau
quotients of \cite[Theorem~2.1]{Hu}. Peters \cite[Example~19]{Pet}
constructed a general type threefold with rational point support; the
vanishing of its Albanese variety and Roitman's torsion theorem give
integral point support. Vial \cite[Section~2.3 and Theorem~2.18]{Via}
proved the integral conclusion for a smooth general type fourfold, and
Laterveer \cite[Theorem~3.1]{Lat} for certain Enriques varieties in
dimensions four and six. Our contribution is the uniform character
criterion, its nonproduct general type realization in every dimension
at least four, and the sharp geometric analysis within this covering
class, rather than a first existence claim in dimension four.

The proofs extend the author's coordinate involution calculation \cite[Lemma~4]{Hu},
within the method of enough automorphisms \cite{IM,Bar},
\cite[Lemma~1.2]{Bau}, by considering all partial covers simultaneously.
When \(d_I\leq n\), the cover retaining the roots indexed by \(I\)
is Fano. Its degree zero Chow group vanishes, so the corresponding
integral character operator annihilates \(A_0(X)\); for
\(I\neq\varnothing\), it annihilates all of
\(\operatorname{CH}_0(X)\). Character summation and transfer yield a
single integer annihilating \(A_0(Y)\), which divisibility over
\(\mathbb{C}\) removes. For products, exterior products of these
operators yield supports formed from selected factors and finite orbits.
Strict transforms transfer the supports to smooth quotient models;
divisibility gives integral surjectivity, and the generalized Mumford
theorem proves optimality. Finally, the vanishing of the Albanese
varieties and Roitman's torsion theorem \cite{Roi80} give unique
divisibility, allowing the natural character-double-cover decomposition
of Theorem~\ref{theorem:character-decomposition}. This is a decomposition
of integral groups, not integral Chow motives. Rational vanishing also
implies integral vanishing in this Albanese-zero setting; these standard
tools are cited, not claimed as new results.

Sections~3 and~4 prove the Chow results, and Section~5 proves
Theorem~\ref{proposition:general type} and describes its examples.
Sections~6--8 develop the geometric refinements, maximal rationally
connected fibrations, and sharp classifications. Automorphisms act by
pushforward on Chow groups and by pullback on forms. Fundamental groups
refer to the associated complex analytic spaces, with base points omitted
when only the isomorphism type matters. Product actions are componentwise
and do not permute factors. All 0-cycle conclusions concern
\(\mathbb{C}\), not universal triviality after arbitrary field extensions.
No argument depends on a general solution of the surface case.\footnote{Entry~040
of the OpenAI mathematics manuscript collection \cite{OAI} reports a proof
that \(A_0(S)\longrightarrow\operatorname{Alb}(S)(\mathbb{C})\) is an
isomorphism for smooth projective surfaces with \(p_g(S)=h^{2,0}(S)=0\).
Here \(\operatorname{Alb}(S)\) is the Albanese variety. Point support
additionally requires \(h^{1,0}(S)=0\). This reported proof has not
been independently verified here and is not used in any argument below.}

\section{Main results}

\subsection{Multiple double covers}

For every positive integer \(r\), write \([r]=\{1,\ldots,r\}\).
Fix integers \(n\geq2\), \(r\geq1\), and \(d_1,\ldots,d_r\geq1\).
For each \(i\in[r]\), choose a nonzero section
\[
F_i\in H^0\bigl(\mathbb{P}^n,
\mathcal{O}_{\mathbb{P}^n}(2d_i)\bigr),
\]
and let \(B_i=(F_i=0)\) be its zero divisor. Throughout the paper, the
covering construction has the following hypotheses:
\begin{equation}\label{equation:branch-hypotheses}
\text{each }B_i\text{ is smooth, and }
B_1+\cdots+B_r\text{ is a reduced simple normal crossings divisor.}
\end{equation}
Thus the branch divisors have no common irreducible component, and the
divisors passing through a given point meet transversally.

Give
\(\mathcal{O}_{\mathbb{P}^n}\oplus
\mathcal{O}_{\mathbb{P}^n}(-d_i)\)
the algebra structure whose multiplication on the second summand is
induced by \(F_i\). Define
\begin{equation}\label{equation:cover-construction}
\begin{split}
\mathcal{A}
&=\bigotimes_{i=1}^r
\bigl(\mathcal{O}_{\mathbb{P}^n}
\oplus\mathcal{O}_{\mathbb{P}^n}(-d_i)\bigr),\\
X&=\operatorname{Spec}_{\mathbb{P}^n}(\mathcal{A}),
\qquad f:X\longrightarrow\mathbb{P}^n.
\end{split}
\end{equation}
The tensor product in \eqref{equation:cover-construction} is a tensor
product of \(\mathcal{O}_{\mathbb{P}^n}\)-algebras. Locally, this construction
is given by the equations \(x_i^2=F_i\). Globally, \(x_i\) is the
tautological section of \(f^*\mathcal{O}_{\mathbb{P}^n}(d_i)\).
Lemma~\ref{lemma:cover-geometry} shows that \(X\) is smooth and connected
and that \(f\) is finite flat of degree \(2^r\).

Let
\[
G=\langle\tau_1,\ldots,\tau_r\rangle
\cong(\mathbb{Z}/2\mathbb{Z})^r,
\]
where \(\tau_i\) changes the sign of \(x_i\) and fixes the other root
coordinates. For a subset \(I\subseteq[r]\), define
\begin{equation}\label{equation:subset-data}
d_I=\sum_{i\in I}d_i,
\qquad
\chi_I(\tau_1^{e_1}\cdots\tau_r^{e_r})
=(-1)^{\sum_{i\in I}e_i},\qquad e_i\in\{0,1\}.
\end{equation}
We use the conventions \(d_{\varnothing}=0\) and
\(\chi_{\varnothing}=1\). For a subgroup \(H\subseteq G\), let
\begin{equation*}
\mathcal{I}(H)=
\{I\subseteq[r]:\chi_I(h)=1\text{ for every }h\in H\}.
\end{equation*}
A smooth projective model of \(X/H\) means a smooth connected projective
variety birational to the finite quotient \(X/H\).

\begin{theorem}\label{theorem:single-cover}
Assume \eqref{equation:branch-hypotheses}. Let \(H\subseteq G\) be any
subgroup, and let \(Y\) be a smooth projective model of \(X/H\). Then
\begin{equation}\label{equation:single-hodge}
H^{p,0}(Y)=0\quad(0<p<n),
\qquad
p_g(Y)=
\sum_{\substack{I\in\mathcal{I}(H)\\d_I\geq n+1}}
\binom{d_I-1}{n}.
\end{equation}
The following conditions are equivalent:
\begin{enumerate}
\item[\textup{(i)}] \(d_I\leq n\) for every \(I\in\mathcal{I}(H)\);
\item[\textup{(ii)}] \(H^{p,0}(Y)=0\) for every \(p>0\);
\item[\textup{(iii)}] the degree map
\(\deg:\operatorname{CH}_0(Y)\longrightarrow\mathbb{Z}\)
is an isomorphism.
\end{enumerate}
In particular, \(p_g(Y)=0\) if and only if
\(\operatorname{CH}_0(Y)\cong\mathbb{Z}\) by the degree map.
\end{theorem}

Only the subsets whose characters are trivial on \(H\) enter the
criterion; there is no bound on the total half-branch degree.
Both divisorial and higher-codimensional fixed loci are allowed.
Since the intermediate holomorphic degrees vanish, the generalized
Mumford theorem gives \(s_0(Y)=n\) when \(p_g(Y)>0\). Together with
the theorem, this proves the full support prediction for these
single-cover quotients.

\subsection{Products}

Fix an integer \(k\geq1\). For each \(a\in\{1,\ldots,k\}\), let
\[
f_a:X_a\longrightarrow\mathbb{P}^{n_a},\qquad n_a\geq2,
\]
be a cover of the form \eqref{equation:cover-construction}, with
\(r_a\) branch divisors of degrees
\(2d_{a,1},\ldots,2d_{a,r_a}\), satisfying
\eqref{equation:branch-hypotheses}. Its group of covering transformations
is denoted by \(G_a\). For \(I\subseteq[r_a]\), let
\[
d_{a,I}=\sum_{i\in I}d_{a,i},
\qquad
\chi_{a,I}:G_a\longrightarrow\{1,-1\}
\]
be the degree sum and character defined as in
\eqref{equation:subset-data}. Set
\begin{equation}\label{equation:admissible-factor}
\mathcal{A}_a=
\{\varnothing\}\cup
\{I\subseteq[r_a]:d_{a,I}\geq n_a+1\}.
\end{equation}
For a tuple \(\boldsymbol{I}=(I_1,\ldots,I_k)\) with
\(I_a\in\mathcal{A}_a\), define
\begin{equation*}
\begin{split}
w(\boldsymbol{I})
&=\sum_{a:I_a\neq\varnothing}n_a,\\
\mu(\boldsymbol{I})
&=\prod_{a:I_a\neq\varnothing}
\binom{d_{a,I_a}-1}{n_a},\\
\chi_{\boldsymbol{I}}(g_1,\ldots,g_k)
&=\prod_{a=1}^k\chi_{a,I_a}(g_a).
\end{split}
\end{equation*}
An empty sum is zero and an empty product is one.
Let
\[
W=X_1\times\cdots\times X_k,
\qquad G=G_1\times\cdots\times G_k,
\]
with the componentwise action of \(G\) on \(W\). For a subgroup
\(H\subseteq G\), define
\begin{equation}\label{equation:product-bound}
\begin{split}
\mathcal{A}(H)
&=\{\boldsymbol{I}\in\mathcal{A}_1\times\cdots\times\mathcal{A}_k:
\chi_{\boldsymbol{I}}|_H=1\},\\
s_H&=\max_{\boldsymbol{I}\in\mathcal{A}(H)}w(\boldsymbol{I}).
\end{split}
\end{equation}
This maximum is defined because the tuple with every entry empty
belongs to \(\mathcal{A}(H)\).

\begin{theorem}\label{theorem:products}
With the notation above, let \(Y\) be a smooth projective model of
\(W/H\). For every \(p\geq0\), one has
\begin{equation}\label{equation:product-hodge}
h^{p,0}(Y)=
\sum_{\substack{\boldsymbol{I}\in\mathcal{A}(H)\\
w(\boldsymbol{I})=p}}\mu(\boldsymbol{I}).
\end{equation}
Moreover, there is a closed reduced algebraic subset \(Z\subseteq Y\)
of dimension at most \(s_H\) such that
\(\operatorname{CH}_0(Z)\longrightarrow\operatorname{CH}_0(Y)\)
is surjective. This dimension is minimal:
\begin{equation}\label{equation:optimal-support}
s_0(Y)=s_H=\max\{p\geq0:H^{p,0}(Y)\neq0\}.
\end{equation}
In particular, the degree map on \(\operatorname{CH}_0(Y)\) is an
isomorphism if and only if \(\mathcal{A}(H)\) consists only of the tuple
with every entry empty.
\end{theorem}

The group \(H\) need not be a product of subgroups of the \(G_a\).
In particular, the theorem applies to diagonal subgroups.

\subsection{Character double covers}

Return to the notation of Section~2.1. For a nonempty subset
\(I\subseteq[r]\), set \(K_I=\ker\chi_I\), and let \(Z_I\)
be a smooth projective model of \(X/K_I\). Equivalently, \(Z_I\)
is a smooth projective model of the double cover of \(\mathbb{P}^n\)
associated with the section
\(\prod_{i\in I}F_i\) of
\(\mathcal{O}_{\mathbb{P}^n}(2d_I)\). This double cover has local
equation
\[
z_I^2=\prod_{i\in I}F_i.
\]
On the double cover, \(z_I\) denotes the tautological section of
the pullback of \(\mathcal{O}_{\mathbb{P}^n}(d_I)\). Define
\begin{equation*}
\mathcal{I}_+(H)=
\{I\in\mathcal{I}(H): I\neq\varnothing,\ d_I\geq n+1\}.
\end{equation*}

\begin{theorem}\label{theorem:character-decomposition}
Let \(Y\) be a smooth projective model of \(X/H\), under the
hypotheses of Theorem~\ref{theorem:single-cover}. For
\(I\in\mathcal{I}_+(H)\), the inclusion \(H\subseteq K_I\)
induces a generically finite rational map
\[
v_I:Y\dashrightarrow Z_I,
\qquad D_I=\deg(v_I)=[K_I:H]=|G/H|/2.
\]
The groups \(A_0(Y)\) and \(A_0(Z_I)\) are uniquely divisible.
Moreover, the homomorphism
\begin{equation}\label{equation:character-decomposition}
\Phi=\bigoplus_{I\in\mathcal{I}_+(H)}v_{I*}:
A_0(Y)\longrightarrow
\bigoplus_{I\in\mathcal{I}_+(H)}A_0(Z_I)
\end{equation}
is an isomorphism, with inverse
\begin{equation*}
\Psi((\alpha_I)_I)=
\sum_{I\in\mathcal{I}_+(H)}D_I^{-1}v_I^*\alpha_I.
\end{equation*}
Here \(v_I^*\) is the pullback on 0-cycles defined by a resolution
of the graph in Section~\ref{subsection:chow-preliminaries}. Division by
\(D_I\) denotes the inverse of multiplication by \(D_I\) on the
uniquely divisible group. The direct sum indexed by the empty set is
the zero group.
\end{theorem}

The isomorphism in \eqref{equation:character-decomposition} concerns
integral Chow groups, without tensoring them with \(\mathbb{Q}\).
It also determines exactly when a further subgroup quotient induces
an isomorphism on 0-cycles; see
Corollary~\ref{corollary:quotient-kernel}.

\subsection{Nonproduct varieties of general type}

The following theorem applies the integral criterion of
Theorem~\ref{theorem:single-cover} to obtain nonproduct varieties of
general type. The geometric assertions refine this Chow theoretic
construction; they are not assumptions in the vanishing criterion.

\begin{theorem}\label{proposition:general type}
For every integer \(n\geq4\), there is a simply connected smooth
projective \(n\)-fold \(Y\) of general type such that
\[
\operatorname{CH}_0(Y)\cong\mathbb{Z}
\qquad\text{and}\qquad
H^{p,0}(Y)=0\quad(p>0).
\]
More precisely, take a cover of \(\mathbb{P}^n\) as in
\eqref{equation:cover-construction} with
\[
r=3,\qquad(d_1,d_2,d_3)=(2,2,n-2),
\]
and let \(\sigma=\tau_1\tau_2\tau_3\).
Then every smooth projective model of
\(Q=X/\langle\sigma\rangle\) has these properties, has canonical
volume four, and is not birational to a product of positive dimensional
smooth projective varieties. In particular, \(\kappa(Y)=n\), and
\(Y\) is not uniruled or rationally connected.
\end{theorem}

The proof is given in Section~5.

\section{Zero-cycles and multiple double covers}

\subsection{Chow theoretic preliminaries}
\label{subsection:chow-preliminaries}

We use the Chow operations and moving lemma of \cite[Chapters~1, 6,
and~16]{Ful}. For a smooth projective variety \(V\), the group
\(A_0(V)\) is divisible \cite[Section~2, Proposition~1]{Roi}.
Consequently, if a subgroup \(B\subseteq\operatorname{CH}_0(V)\)
contains a degree-one class and
\begin{equation*}
N\operatorname{CH}_0(V)\subseteq B
\end{equation*}
for one positive integer \(N\), then \(B=\operatorname{CH}_0(V)\):
the quotient is divisible and annihilated by \(N\). We shall use this
bounded exponent argument repeatedly. Also, \(s_0(V)=0\) if and only
if the degree map is an isomorphism, since a finitely generated divisible
degree zero group is zero. If \(H^{1,0}(V)=0\), then \(A_0(V)\) is
uniquely divisible by Roitman's torsion theorem \cite{Roi80}.

For a dominant rational map \(u:V\dashrightarrow Y\) between smooth
projective varieties, choose a resolution of its graph
\[
V\xleftarrow{\ b\ }\widetilde V\xrightarrow{\ v\ }Y.
\]
Birational invariance of 0-cycles \cite[Example~16.1.11]{Ful}
allows us to define
\begin{equation}\label{equation:rational-pushforward}
u_*=v_*b_*^{-1}:\operatorname{CH}_0(V)\longrightarrow
\operatorname{CH}_0(Y).
\end{equation}
This map preserves degree and is surjective over \(\mathbb{C}\),
since every closed point of \(Y\) lifts through \(v\). For a
generically finite \(u\) of degree \(d\), use the usual pullback
between smooth Chow rings in codimension \(\dim Y\) to define
\begin{equation*}
u^*=b_*v^*:\operatorname{CH}_0(Y)\longrightarrow
\operatorname{CH}_0(V).
\end{equation*}
These are the graph-correspondence operations. They are independent of
the resolution and functorial for dominant generically finite rational
maps. The projection and transfer formulas give
\begin{equation*}
u_*u^*=d\operatorname{id}.
\end{equation*}
If \(u\) is a rational quotient by an effective finite group \(H\),
followed by a choice of smooth projective model, then
\begin{equation*}
u^*u_*=\sum_{h\in H}h_*,\qquad u_*h_*=u_*.
\end{equation*}
We use these standard identities as in \cite[Lemma~3]{Hu} and
\cite[Chapters~6 and~16]{Ful}. In the generically finite case, the
moving lemma places cycles in a common finite \'etale locus, where
these identities and functoriality are the usual ones for finite covers.
No flatness of \(v\) is assumed. For any automorphism \(g\) with
\(u\circ g=u\), the graph construction gives \(u_*g_*=u_*\).

\begin{lemma}\label{lemma:transport-support}
Use the graph resolution in Section~\ref{subsection:chow-preliminaries}, and let \(U\subseteq V\)
be a dense open subset over which \(b\) is an isomorphism. Let
\(T\subseteq V\) be a nonempty closed reduced algebraic subset such
that every irreducible component of \(T\) meets \(U\). Denote its
inclusion by \(j:T\hookrightarrow V\). Suppose that, for some positive
integer \(N\),
\begin{equation}\label{equation:transport-assumption}
N\operatorname{CH}_0(Y)
\subseteq u_*j_*\operatorname{CH}_0(T).
\end{equation}
Then \(\operatorname{CH}_0(Y)\) is supported on a closed reduced
algebraic subset of dimension at most \(\dim T\).
\end{lemma}

\begin{proof}
Let \(T'\subseteq\widetilde{V}\) be the union, with its reduced
structure, of the strict transforms of the irreducible components of
\(T\). The restriction \(b_T:T'\longrightarrow T\) is proper and
surjective, and \(\dim T'=\dim T\). Pushforward
\(\operatorname{CH}_0(T')\longrightarrow\operatorname{CH}_0(T)\)
is surjective, since every closed point has a closed point above it.
As \(b_*\) is an isomorphism, the identity
\(b_*j'_* = j_*b_{T*}\), where
\(j':T'\hookrightarrow\widetilde{V}\), gives
\[
b_*^{-1}j_*\operatorname{CH}_0(T)
\subseteq j'_*\operatorname{CH}_0(T').
\]
Set \(Z=v(T')\) with its reduced structure. It is closed, nonempty,
and of dimension at most \(\dim T\). By
\eqref{equation:rational-pushforward} and
\eqref{equation:transport-assumption}, the image of
\(\operatorname{CH}_0(Z)\) in \(\operatorname{CH}_0(Y)\) contains
\(N\operatorname{CH}_0(Y)\). This image contains a degree-one class.
The bounded exponent argument of Section~\ref{subsection:chow-preliminaries} proves the assertion.
\end{proof}

For holomorphic forms on smooth quotient models, we shall also use the
standard formula
\begin{equation}\label{equation:quotient-forms}
h^{p,0}(Y)=\dim_{\mathbb{C}}H^{p,0}(V)^H\qquad(p\geq0),
\end{equation}
where a finite group \(H\) acts on a smooth projective variety \(V\)
and \(Y\) is any smooth projective model of \(V/H\). Finite quotients
have rational singularities by \cite{Bou}. Thus resolution preserves
\(H^p(\mathcal{O})\), and the finite quotient identifies it with
\(H^p(V,\mathcal{O}_V)^H\). Formula~\eqref{equation:quotient-forms}
is the corresponding statement under Hodge symmetry.

\subsection{The covering algebra and intermediate covers}

For \(I\subseteq[r]\), let \(X_I\) be the partial cover obtained by
retaining precisely the roots indexed by \(I\):
\begin{equation*}
X_I=\operatorname{Spec}_{\mathbb{P}^n}
\left(\bigotimes_{i\in I}
\bigl(\mathcal{O}_{\mathbb{P}^n}
\oplus\mathcal{O}_{\mathbb{P}^n}(-d_i)\bigr)\right).
\end{equation*}
Write \(f_I:X_I\longrightarrow\mathbb{P}^n\) for its structure map
and \(q_I:X\longrightarrow X_I\) for the natural morphism. Thus
\(X_{[r]}=X\), while \(X_{\varnothing}=\mathbb{P}^n\).
Set \(I^c=[r]\setminus I\) and
\(G_{I^c}=\langle\tau_j:j\in I^c\rangle\).

\begin{lemma}\label{lemma:cover-geometry}
Under \eqref{equation:branch-hypotheses}, every \(X_I\) is smooth,
connected, and projective. The morphism \(f_I\) is finite flat of
degree \(2^{|I|}\), and \(q_I\) is the finite flat quotient by
\(G_{I^c}\). Moreover,
\begin{equation}\label{equation:cover-sheaf}
f_*\mathcal{O}_X
=\bigoplus_{I\subseteq[r]}
\mathcal{O}_{\mathbb{P}^n}(-d_I),
\end{equation}
where the summand indexed by \(I\) has character \(\chi_I\), and
\begin{equation}\label{equation:canonical-cover}
\omega_{X_I}\cong
f_I^*\mathcal{O}_{\mathbb{P}^n}(d_I-n-1).
\end{equation}
In addition, \(\pi_1(X_I)=1\); if \(n\geq3\), then
\(\operatorname{Pic}(X_I)\cong\mathbb{Z}\).
In particular, \(d_I\leq n\) implies
\(\operatorname{CH}_0(X_I)\cong\mathbb{Z}\) by the degree map.
\end{lemma}

\begin{proof}
The square-free root monomials form a local basis of each defining
algebra. They give finite flatness, the degrees, the eigensheaf formula
\eqref{equation:cover-sheaf}, and the description of \(q_I\) as a
quotient. At a point of the base, the vanishing branch equations are
part of a regular system of parameters. Adjoining their square roots
replaces those parameters by root coordinates; adjoining square roots
of units is \'etale. Hence every partial cover is smooth. The
negative-degree summands of its covering algebra have no global
sections, so \(H^0(X_I,\mathcal{O}_{X_I})=\mathbb{C}\); thus it is
connected and irreducible. It is projective because it is finite over
\(\mathbb{P}^n\).

The canonical bundle formula \eqref{equation:canonical-cover} is the
standard ramification formula for abelian covers \cite{Par}, applied
to the branch line bundles \(\mathcal{O}_{\mathbb{P}^n}(d_i)\).
For the topological assertions, realize \(X_I\) as
\[
x_i^2=F_i(y_0,\ldots,y_n)\quad(i\in I)
\quad\text{in}\quad
\mathbb{P}\bigl(\underbrace{1,\ldots,1}_{n+1},(d_i)_{i\in I}\bigr).
\]
The equations are a regular sequence, being monic in distinct root
variables. They imply that at every projective point some weight-one
coordinate \(y_j\) is nonzero. In particular, \(X_I\) is contained
in the smooth locus of the weighted ambient space. On that chart the punctured affine cone is the
product of the corresponding smooth chart of \(X_I\) with
\(\mathbb{G}_m\), the multiplicative group. Consequently the weighted
complete intersection is quasismooth. The assertions on the Picard
group and the topological fundamental group follow directly from
\cite[Theorem~3.2.4(i), (ii)$'$]{Dol}. For \(I=\varnothing\), use
\(X_I=\mathbb{P}^n\).

If \(d_I\leq n\), the displayed canonical bundle formula makes
\(X_I\) Fano. Its rational connectedness follows from
\cite[Corollaire~3.2]{Cam}, \cite{KMM}, and hence its degree map on
0-cycles is an isomorphism; see \cite{AK}.
\end{proof}

The isomorphism \eqref{equation:canonical-cover} is an isomorphism of
line bundles. When computing the group action on holomorphic forms,
we use the eigensheaf decomposition
\eqref{equation:cover-sheaf}, so that no choice of linearization of the
right-hand side of \eqref{equation:canonical-cover} is involved.

\begin{lemma}\label{lemma:character-operators}
For each \(I\subseteq[r]\), define an element of the integral group
ring \(\mathbb{Z}[G]\) by
\begin{equation}\label{equation:character-numerator}
\Theta_I=\sum_{g\in G}\chi_I(g)g
=\prod_{i\in I}(1-\tau_i)
\prod_{j\in I^c}(1+\tau_j).
\end{equation}
An element of \(\mathbb{Z}[G]\) acts on
\(\operatorname{CH}_0(X)\) by the corresponding sum of proper
pushforwards. Then:
\begin{enumerate}
\item[\textup{(i)}] if \(d_I\leq n\), then
\(\Theta_I A_0(X)=0\);
\item[\textup{(ii)}] if \(I\neq\varnothing\) and \(d_I\leq n\),
then \(\Theta_I\operatorname{CH}_0(X)=0\);
\item[\textup{(iii)}] for any points \(x,a\in X\),
\(\Theta_{\varnothing}[x]=\Theta_{\varnothing}[a]\);
\item[\textup{(iv)}] one has
\(g\Theta_I=\chi_I(g)\Theta_I\) for \(g\in G\), and
\(\sum_{I\subseteq[r]}\Theta_I=|G|\cdot1\) in \(\mathbb{Z}[G]\).
\end{enumerate}
\end{lemma}

\begin{proof}
Finite flat pullback along \(q_I\) is defined, and
\begin{equation}\label{equation:integral-transfer}
q_I^*q_{I*}=\sum_{g\in G_{I^c}}g_*
=\prod_{j\in I^c}(\operatorname{id}+\tau_{j*})
\quad\text{on }\operatorname{CH}_0(X).
\end{equation}
This is the finite-flat quotient transfer formula \cite{Ful}.
It can also be checked on the free locus and extended to
\(\operatorname{CH}_0(X)\) by the moving lemma; this avoids any
assumption that the quotient is unramified everywhere.

When \(d_I\leq n\), the group \(A_0(X_I)\) vanishes by
Lemma~\ref{lemma:cover-geometry}. Since \(q_{I*}\) preserves degree,
\eqref{equation:integral-transfer} kills \(A_0(X)\). Factoring
\(\Theta_I\) as in \eqref{equation:character-numerator} proves (i).

To prove (ii), let \(\overline{\tau}_i\) denote the involution induced
by \(\tau_i\) on \(X_I\), for \(i\in I\). Equivariance of flat
pullback gives
\[
\Theta_I=q_I^*
\left(\prod_{i\in I}
(\operatorname{id}-\overline{\tau}_{i*})\right)q_{I*}
\quad\text{on }\operatorname{CH}_0(X).
\]
Every \(\overline{\tau}_{i*}\) is the identity on
\(\operatorname{CH}_0(X_I)\cong\mathbb{Z}\), because it preserves
degree. The middle product is therefore zero when \(I\neq\varnothing\).

For (iii), use
\(\Theta_{\varnothing}=f^*f_*\) and
\([f(x)]=[f(a)]\) in
\(\operatorname{CH}_0(\mathbb{P}^n)\).
Finally, the first identity in (iv) follows by reindexing the character
sum, and the second follows from character orthogonality. All these
identities hold with integral coefficients.
\end{proof}

\begin{proof}[Proof of Theorem~\ref{theorem:single-cover}]

The cohomology of line bundles on projective space and
\eqref{equation:cover-sheaf} give
\begin{equation*}
H^p(X,\mathcal{O}_X)=0\quad(0<p<n),
\qquad
H^n(X,\mathcal{O}_X)
=\bigoplus_{I\subseteq[r]}
H^n\bigl(\mathbb{P}^n,
\mathcal{O}_{\mathbb{P}^n}(-d_I)\bigr).
\end{equation*}
The summand indexed by \(I\) has character \(\chi_I\). Serre duality
on \(\mathbb{P}^n\) shows that its dimension is
\(\binom{d_I-1}{n}\) if \(d_I\geq n+1\), and zero otherwise.
Complex conjugation identifies the representations on
\(H^{0,n}(X)\) and \(H^{n,0}(X)\) after conjugating the eigenvalues.
Since the characters take values in \(\{1,-1\}\), the character
labels do not change. Formula~\eqref{equation:quotient-forms} now proves
\eqref{equation:single-hodge}, and hence the equivalence of (i) and (ii).

Assume (i), and set
\(N_H=\sum_{h\in H}h\in\mathbb{Z}[G]\).
The characters \(\chi_I\) with \(I\in\mathcal{I}(H)\) are exactly
the characters of \(G/H\). Consequently,
\begin{equation}\label{equation:subgroup-orthogonality}
\sum_{I\in\mathcal{I}(H)}\Theta_I
=|G/H|N_H
\qquad\text{in }\mathbb{Z}[G].
\end{equation}
Indeed, the coefficient of \(g\in G\) on the left is \(|G/H|\)
when \(g\in H\), and zero otherwise. By
Lemma~\ref{lemma:character-operators}(i), the left-hand side kills
\(A_0(X)\). Let \(u:X\dashrightarrow Y\) be the rational map
induced by the quotient. Section~\ref{subsection:chow-preliminaries}
gives a surjection \(A_0(X)\longrightarrow A_0(Y)\), and
\(u_*N_H=|H|u_*\). Pushing forward
\eqref{equation:subgroup-orthogonality} therefore gives
\[
|G|A_0(Y)=0.
\]
Divisibility implies
\(A_0(Y)=0\). This proves (iii).

Finally, (iii) says that \(\operatorname{CH}_0(Y)\) is supported on a
point. The generalized Mumford theorem
\cite[Theorem~3.13]{Voi} gives (ii), completing the proof.
\end{proof}

\subsection{Transfers and character double covers}

\begin{proof}[Proof of Theorem~\ref{theorem:character-decomposition}]
By Theorem~\ref{theorem:single-cover}, the varieties \(X\), \(Y\),
and \(Z_I\) have \(H^{1,0}=0\), since \(n\geq2\). Their degree zero
Chow groups are therefore uniquely divisible, as recalled in
Section~\ref{subsection:chow-preliminaries}. We use their canonical
\(\mathbb{Q}\)-vector space structures; every group homomorphism
between them respects these structures.

Write \(A=A_0(X)\). For a subgroup \(L\subseteq G\), write
\(A^L\) for its subgroup of \(L\)-invariant elements, and define
\(N_L=\sum_{g\in L}g_*\). Let \(u_H:X\dashrightarrow Y\)
be the rational quotient map. The transfer identities imply
\[
u_{H*}u_H^*=|H|\operatorname{id},\qquad
u_H^*u_{H*}=N_H.
\]
Since \(u_{H*}\) is surjective, the norm identity shows that the
image of \(u_H^*\) is contained in \(A^H\).
Consequently,
\begin{equation}\label{equation:invariant-chow-group}
u_H^*:A_0(Y)\xrightarrow{\ \cong\ }A^H,
\end{equation}
with inverse \(|H|^{-1}u_{H*}\). Indeed, the first identity proves
injectivity, and for \(a\in A^H\) the second gives
\(u_H^*(|H|^{-1}u_{H*}a)=a\).

Define \(e_I=|G|^{-1}\Theta_I\) as an endomorphism of \(A\), where
\(\Theta_I\) is the integral character operator in
\eqref{equation:character-numerator}. Unique divisibility makes this
definition meaningful on the integral group \(A\). Character
orthogonality gives pairwise orthogonal idempotents with sum the
identity. Lemma~\ref{lemma:character-operators} gives
\(e_IA=0\) whenever \(d_I\leq n\), including
\(I=\varnothing\). It follows that
\begin{equation}\label{equation:invariant-character-splitting}
A^H=\bigoplus_{I\in\mathcal{I}_+(H)}e_IA.
\end{equation}
For nonempty \(I\), the only characters trivial on \(K_I\) are
\(1\) and \(\chi_I\). Since the trivial summand of \(A\) is
zero, the rational quotient map \(u_I:X\dashrightarrow Z_I\)
therefore gives
\begin{equation}\label{equation:double-cover-summand}
u_I^*:A_0(Z_I)\xrightarrow{\ \cong\ }A^{K_I}=e_IA.
\end{equation}

For \(I\in\mathcal{I}_+(H)\), the equality
\(u_I=v_I\circ u_H\) and Section~\ref{subsection:chow-preliminaries}
give
\(u_I^*=u_H^*v_I^*\) and
\(u_{I*}=v_{I*}u_{H*}\). If \(y\in A_0(Y)\), then
\begin{align*}
u_H^*\bigl(D_I^{-1}v_I^*v_{I*}y\bigr)
&=(D_I|H|)^{-1}u_I^*u_{I*}u_H^*y\\
&=|K_I|^{-1}N_{K_I}u_H^*y
=e_Iu_H^*y.
\end{align*}
The last equality holds because the trivial character summand of
\(A\) vanishes. Summing over \(I\) and using
\eqref{equation:invariant-character-splitting} proves
\(\Psi\Phi=\operatorname{id}\). On the other hand, under
\eqref{equation:invariant-chow-group}, the map \(\Psi\) is the
direct sum of the isomorphisms
\(D_I^{-1}u_I^*\) in \eqref{equation:double-cover-summand}.
Thus \(\Psi\) is an isomorphism, and \(\Phi=\Psi^{-1}\).
\end{proof}

\begin{corollary}\label{corollary:quotient-kernel}
Let \(H\subseteq K\subseteq G\), let \(Y_H\) and \(Y_K\)
be smooth projective models of \(X/H\) and \(X/K\), and let
\(w:Y_H\dashrightarrow Y_K\) be the induced rational map.
Then
\[
\ker\bigl(w_*:A_0(Y_H)\longrightarrow A_0(Y_K)\bigr)
\cong
\bigoplus_{I\in\mathcal{I}_+(H)\setminus\mathcal{I}_+(K)}
A_0(Z_I).
\]
The map on degree zero groups is split by
\((\deg w)^{-1}w^*\). Moreover,
\[
w_*:\operatorname{CH}_0(Y_H)\longrightarrow
\operatorname{CH}_0(Y_K)\text{ is an isomorphism}
\quad\Longleftrightarrow\quad
\mathcal{I}_+(H)=\mathcal{I}_+(K).
\]
\end{corollary}

\begin{proof}
For \(I\in\mathcal{I}_+(K)\), the quotient map from \(Y_H\)
to \(Z_I\) factors through \(w\). Functoriality therefore
identifies \(w_*\), under Theorem~\ref{theorem:character-decomposition},
with the projection retaining these summands. The transfer degree
identity gives the asserted splitting. Every summand indexed by
\(\mathcal{I}_+(H)\) is nonzero: by
\eqref{equation:single-hodge},
\(p_g(Z_I)=\binom{d_I-1}{n}>0\), and the generalized Mumford
theorem excludes \(A_0(Z_I)=0\). Thus the kernel vanishes exactly
when the two index sets agree. Finally, \(w_*\) preserves degree
and is surjective, so the same criterion applies to the full Chow
groups.
\end{proof}

\section{Products and the dimension of a support}

\begin{proof}[Proof of Theorem~\ref{theorem:products}]

For the factor \(X_a\), let \(\Theta_{a,I}\) denote the character operator
of Lemma~\ref{lemma:character-operators}. The decomposition of
holomorphic forms on \(X_a\) has a one-dimensional trivial summand
in degree zero and, in degree \(n_a\), a summand of dimension
\(\binom{d_{a,I}-1}{n_a}\) with character \(\chi_{a,I}\) for each
nonempty \(I\in\mathcal{A}_a\). There are no holomorphic forms in
the intervening degrees. The K\"unneth formula and
\eqref{equation:quotient-forms} prove
\eqref{equation:product-hodge}.

We now construct a support. Let \(u:W\dashrightarrow Y\) be the
rational quotient map, and choose a resolution of its graph as in
Section~\ref{subsection:chow-preliminaries}. Let \(U\subseteq W\) be a
dense open subset over which the birational morphism from that
resolution is an isomorphism. Choose a point
\(o=(o_1,\ldots,o_k)\in W\) such that its entire \(G\)-orbit is
contained in \(U\). Such a point exists, since
\(\bigcap_{g\in G}g^{-1}U\) is a dense open subset of \(W\).
Define
\begin{equation*}
\beta_a=\Theta_{a,\varnothing}[o_a]
=\sum_{g\in G_a}[g(o_a)],
\qquad e_a=|G_a|,
\qquad M=\prod_{a=1}^k e_a.
\end{equation*}
The cycle \(\beta_a\) is supported on the finite orbit
\(G_a\cdot o_a\) and has degree \(e_a\).

For a point \(x=(x_1,\ldots,x_k)\in W\), the identity
\(\sum_I \Theta_{a,I}=e_a\operatorname{id}\) in each factor gives
\begin{equation}\label{equation:product-expansion}
M[x]=\sum_{\boldsymbol{I}\in
\mathcal{A}_1\times\cdots\times\mathcal{A}_k}
z_{\boldsymbol{I}}(x)
\quad\text{in }\operatorname{CH}_0(W),
\end{equation}
where
\[
z_{\boldsymbol{I}}(x)
=\Theta_{1,I_1}[x_1]\boxtimes\cdots\boxtimes \Theta_{k,I_k}[x_k].
\]
The symbol \(\boxtimes\) denotes the exterior product of 0-cycles,
with factors in the indicated order. All nonempty subsets outside
\(\mathcal{A}_a\) give zero by
Lemma~\ref{lemma:character-operators}(ii), which explains the indexing
in \eqref{equation:product-expansion}. If \(I_a=\varnothing\), the
corresponding factor in \(z_{\boldsymbol{I}}(x)\) is \(\beta_a\),
by part (iii) of that lemma.

Let \(N_H=\sum_{h\in H}h_*\) act on
\(\operatorname{CH}_0(W)\). Each term in
\eqref{equation:product-expansion} transforms by
\(\chi_{\boldsymbol{I}}\). The sum of a nontrivial character on a
finite group is zero, so
\begin{equation}\label{equation:norm-expansion}
M N_H[x]
=|H|\sum_{\boldsymbol{I}\in\mathcal{A}(H)}
z_{\boldsymbol{I}}(x).
\end{equation}
Define a closed reduced algebraic subset
\begin{equation*}
T=\bigcup_{\boldsymbol{I}\in\mathcal{A}(H)}
\prod_{a=1}^k T_{a,I_a},
\qquad
T_{a,I_a}=
\begin{cases}
X_a,& I_a\neq\varnothing,\\
G_a\cdot o_a,& I_a=\varnothing.
\end{cases}
\end{equation*}
Then \(\dim T=s_H\), and
\eqref{equation:norm-expansion} shows that
\(M N_H\operatorname{CH}_0(W)\) lies in the image of
\(\operatorname{CH}_0(T)\). Every irreducible component of \(T\)
contains a point of the chosen \(G\)-orbit, and therefore meets
\(U\).

By Section~\ref{subsection:chow-preliminaries}, the map \(u_*\) is
surjective and satisfies \(u_*N_H=|H|u_*\). We obtain
\[
M|H|\operatorname{CH}_0(Y)
\subseteq u_*\operatorname{im}
\bigl(\operatorname{CH}_0(T)\longrightarrow
\operatorname{CH}_0(W)\bigr).
\]
Lemma~\ref{lemma:transport-support} produces a closed subset of
\(Y\) of dimension at most \(s_H\) supporting
\(\operatorname{CH}_0(Y)\). This proves the upper bound with
integral coefficients.

Every integer \(\mu(\boldsymbol{I})\) in
\eqref{equation:product-hodge} is positive. Hence
\(H^{s_H,0}(Y)\neq0\), while \(H^{p,0}(Y)=0\) for \(p>s_H\).
If \(s_H>0\), the generalized Mumford theorem
\cite[Theorem~3.13]{Voi} rules out a support of dimension strictly
less than \(s_H\). For \(s_H=0\), the lower bound is automatic.
This proves \eqref{equation:optimal-support}. Since all \(n_a\)
are positive, \(s_H=0\) is equivalent to the last condition in the
theorem. The criterion for zero-dimensional support in Section~\ref{subsection:chow-preliminaries} completes the proof.
\end{proof}

\subsection{Finite groups acting through sign characters}

The preceding support argument has a useful form that does not require
the varieties to be multiple double covers or the acting groups to be
abelian.

\begin{proposition}\label{proposition:scalar-products}
For \(i\in\{1,\ldots,k\}\), let \(V_i\) be a smooth connected
projective variety of dimension \(m_i>0\), let \(\Gamma_i\) be a finite
group acting on \(V_i\), and let
\(\epsilon_i:\Gamma_i\longrightarrow\{1,-1\}\) be a character.
Assume that
\begin{equation*}
g_*|_{A_0(V_i)}=\epsilon_i(g)\operatorname{id}
\qquad\text{for every }g\in\Gamma_i.
\end{equation*}
Let \(L\subseteq\Gamma_1\times\cdots\times\Gamma_k\) be a subgroup,
acting componentwise, and let \(Y\) be a smooth projective model of
\((V_1\times\cdots\times V_k)/L\). For \(J\subseteq[k]\), define
\[
m_J=\sum_{i\in J}m_i,
\qquad
\epsilon_J(g_1,\ldots,g_k)=\prod_{i\in J}\epsilon_i(g_i),
\qquad
t_L=\max_{\epsilon_J|_L=1}m_J.
\]
Then \(s_0(Y)\leq t_L\).

Suppose in addition that, for every \(i\),
\begin{equation}\label{equation:scalar-hodge-hypotheses}
\begin{gathered}
H^{p,0}(V_i)=0\quad(0<p<m_i),\qquad
b_i=h^{m_i,0}(V_i)>0,\\
g^*|_{H^{m_i,0}(V_i)}=\epsilon_i(g)\operatorname{id}
\quad(g\in\Gamma_i).
\end{gathered}
\end{equation}
Then
\begin{equation}\label{equation:scalar-product-hodge}
h^{p,0}(Y)=
\sum_{\substack{J\subseteq[k]\\m_J=p,\ \epsilon_J|_L=1}}
\prod_{i\in J}b_i,
\qquad s_0(Y)=t_L.
\end{equation}
Under these additional hypotheses, the degree map on
\(\operatorname{CH}_0(Y)\) is an isomorphism if and only if the map
\begin{equation}\label{equation:combined-sign}
L\longrightarrow\{1,-1\}^k,
\qquad(g_1,\ldots,g_k)\longmapsto
(\epsilon_1(g_1),\ldots,\epsilon_k(g_k))
\end{equation}
is surjective.
\end{proposition}

\begin{proof}
Write \(V=\prod_iV_i\) and \(\Gamma=\prod_i\Gamma_i\). Choose
\(o=(o_i)_i\) so that its \(\Gamma\)-orbit is contained in the
open subset of \(V\) over which the birational morphism in a graph
resolution of \(V\dashrightarrow Y\) is an isomorphism.
For \(e_i=|\Gamma_i|\), define
\[
\beta_i=\sum_{g\in\Gamma_i}[g(o_i)],\qquad
\gamma_i(x_i)=e_i[x_i]-\beta_i\in A_0(V_i).
\]
Then \(\beta_i\) is invariant and \(\gamma_i(x_i)\) has character
\(\epsilon_i\). Expanding the exterior product gives
\begin{equation*}
\left(\prod_i e_i\right)[x]
=\boxtimes_{i=1}^k\bigl(\gamma_i(x_i)+\beta_i\bigr).
\end{equation*}
The term using \(\gamma_i\) exactly for \(i\in J\) has character
\(\epsilon_J\). Averaging over \(L\) eliminates the noninvariant
terms. Each surviving term is supported on a product of the factors
\(V_i\), for \(i\in J\), and the finite orbits \(\Gamma_i o_i\)
otherwise. This support has dimension at most \(t_L\), and each
component meets the chosen open. The proof of
Theorem~\ref{theorem:products}, using
Lemma~\ref{lemma:transport-support} with the integer
\(|L|\prod_i e_i\), therefore applies without change.

Under \eqref{equation:scalar-hodge-hypotheses}, the K\"unneth formula
and \eqref{equation:quotient-forms} give
\eqref{equation:scalar-product-hodge}. The generalized Mumford theorem
\cite[Theorem~3.13]{Voi} gives equality in the support bound. Finally,
\(t_L=0\) precisely when no nonempty subset character is trivial on
\(L\). By character duality for \(\{1,-1\}^k\), this is equivalent
to surjectivity in \eqref{equation:combined-sign}.
\end{proof}

\begin{remark}\label{remark:scope}
The actions are componentwise. An action permuting factors would require
additional relations between exterior products of 0-cycles; the scalar
action hypothesis does not supply such relations.
\end{remark}

\section{Sign characters and examples}

\subsection{A single character on degree zero 0-cycles}

Retain the notation of Section~2.1, and write
\[
\delta=d_1+\cdots+d_r,
\qquad d_{\min}=\min_{1\leq i\leq r}d_i,
\qquad\chi=\chi_{[r]}.
\]
We call \(\chi\) the total sign character. An element \(g\in G\)
has \(\chi(g)=-1\) precisely when it changes an odd number of root
coordinates.

\begin{corollary}\label{corollary:parity}
Assume \eqref{equation:branch-hypotheses} and
\begin{equation*}
\delta\leq n+d_{\min}.
\end{equation*}
Then
\begin{equation}\label{equation:parity-action}
g_*|_{A_0(X)}=\chi(g)\operatorname{id}
\qquad(g\in G).
\end{equation}
If \(\delta\leq n\), the degree map is an isomorphism for every
smooth projective model of every subgroup quotient \(X/H\).
If instead \(n+1\leq\delta\leq n+d_{\min}\), then
\[
h^{n,0}(X)=\binom{\delta-1}{n}>0,
\]
the action on \(H^{n,0}(X)\) is given by \(\chi\), and for a
smooth projective model \(Y\) of \(X/H\) the following are
equivalent:
\begin{enumerate}
\item[\textup{(i)}] \(\chi|_H\neq1\);
\item[\textup{(ii)}] \(H^{p,0}(Y)=0\) for every \(p>0\);
\item[\textup{(iii)}] the degree map on \(\operatorname{CH}_0(Y)\)
is an isomorphism.
\end{enumerate}
\end{corollary}

\begin{proof}
Every proper subset \(I\subsetneq[r]\) satisfies
\(d_I\leq\delta-d_{\min}\leq n\). On the uniquely divisible group
\(A_0(X)\), the character decomposition in the proof of
Theorem~\ref{theorem:character-decomposition} therefore has only the
possible summand indexed by \([r]\). Its character is \(\chi\),
which proves \eqref{equation:parity-action}. If \(\delta\leq n\),
this last summand also vanishes, and
Theorem~\ref{theorem:single-cover} applies to every subgroup. Otherwise,
\([r]\) is also the only subset contributing a positive-degree
holomorphic form, with multiplicity \(\binom{\delta-1}{n}\).
The same theorem gives the stated equivalences.
\end{proof}

\Needspace{5\baselineskip}
We first recall the specified Calabi--Yau quotients studied in \cite{Hu}.

\begin{example}\label{example:hu}
Take \(r=n+1\) and \(d_i=1\). Then \(X\) is the smooth complete
intersection
\[
x_i^2=F_i(y_0,\ldots,y_n)\quad(1\leq i\leq n+1)
\qquad\text{in }\mathbb{P}^{2n+1}.
\]
It has \(\omega_X\cong\mathcal{O}_X\),
\(H^{p,0}(X)=0\) for \(0<p<n\), and \(h^{n,0}(X)=1\).
For even \(n\), the involution \(\tau_1\cdots\tau_{n+1}\) is
free. For odd \(n\geq3\), the involution \(\tau_1\cdots\tau_n\)
has \(2^{n+1}\) fixed points, with local quotient type
\(\frac{1}{2}(1,\ldots,1)\). The degree map on every smooth projective
model of either specified quotient is an isomorphism. For general
defining quadrics this is \cite[Theorem~2.1 and Lemmas~5--6]{Hu}.
Under the branch hypotheses of the present paper, it follows from
Corollary~\ref{corollary:parity}, which also treats every subgroup with
nontrivial total sign character. Thus the cited generic construction
and the present branch-uniform statement are distinguished.

For both specified quotients, the canonical and fundamental group
formulas of Section~\ref{section:birational-type} give
\(\kappa(Y)=0\) and \(P_2(Y)=1\). Their fundamental groups are
\(\mathbb{Z}/2\mathbb{Z}\) in even dimension and trivial in odd
dimension. In even dimension the quotient is smooth and its canonical
line bundle has order two. In odd dimension the quotient has terminal
singularities. On the resolution obtained by blowing up the fixed points
upstairs and then taking the quotient, one has
\begin{equation*}
K_Y\sim_{\mathbb{Q}}\rho^*K_Q+
\frac{n-2}{2}\sum_{j=1}^{2^{n+1}}D_j
\sim_{\mathbb{Q}}\frac{n-2}{2}\sum_{j=1}^{2^{n+1}}D_j,
\end{equation*}
where \(\rho:Y\longrightarrow Q\) is the resolution and
\(D_j\cong\mathbb{P}^{n-1}\) are its exceptional divisors.
This is the blow-up formula for the local models of
\cite[Remark~3]{Hu}. In particular, the canonical class of this smooth
resolution is not torsion, although \(K_Q\sim_{\mathbb{Q}}0\).
The models are not uniruled. In contrast, quotienting \(X\) by any
individual \(\tau_i\) gives a smooth Fano partial cover.
\end{example}

\Needspace{5\baselineskip}
The next example is a diagonal quotient of three K3 surfaces with positive dimensional minimal support.

\begin{example}\label{example:three-k3}
For \(i=1,2,3\), let \(S_i\) be a cover as in
Example~\ref{example:hu} with \(n=2\), and choose a covering
involution \(\sigma_i\) with total sign \(-1\).
Each \(S_i\) is a smooth surface with trivial canonical bundle and
\(H^1(S_i,\mathcal{O}_{S_i})=0\), hence a K3 surface. The involution
\(\sigma_i\) acts as minus the identity on \(A_0(S_i)\) and on
\(H^{2,0}(S_i)\). Let \(Y\) be a smooth projective model of
\[
(S_1\times S_2\times S_3)/
\langle(\sigma_1,\sigma_2,\sigma_3)\rangle.
\]
The invariant subset characters correspond to the empty subset and
the three subsets of cardinality two. Therefore
\[
h^{4,0}(Y)=3,
\qquad
H^{p,0}(Y)=0\quad(p>0,\ p\neq4),
\qquad s_0(Y)=4.
\]
In particular, \(\operatorname{CH}_0(Y)\) is supported on a
four-dimensional closed subset, and it is not supported on a subset
of dimension at most three.
No individual coordinate involution belongs to the diagonal subgroup.
Thus Proposition~\ref{proposition:product-kodaira} gives \(\kappa(Y)=0\).
These sixfolds are not uniruled, and their 0-cycle group is not
isomorphic to \(\mathbb{Z}\).
\end{example}

\subsection{Quotients of general type}

\begin{proof}[Proof of Theorem~\ref{proposition:general type}]
Take general branch hypersurfaces of degrees \(4,4,2n-4\).
They satisfy \eqref{equation:branch-hypotheses} by Bertini's theorem. Then
\(\delta=n+2\), \(d_{\min}=2\), and
\[
\omega_X\cong f^*\mathcal{O}_{\mathbb{P}^n}(1),\qquad
h^{n,0}(X)=n+1.
\]
The total sign of \(\sigma\) is \(-1\), so
Corollary~\ref{corollary:parity} gives the 0-cycle and Hodge assertions.
No individual coordinate involution belongs to \(\langle\sigma\rangle\).
Proposition~\ref{proposition:single-kodaira} therefore gives
\(K_Q\sim_{\mathbb{Q}}\overline f^{\,*}A\), which is ample, and
shows that every smooth projective model is of general type. Its
canonical volume is \(2^3/2=4\) by
Proposition~\ref{proposition:plurigenera}, and it has no birational
product decomposition by Proposition~\ref{proposition:nonproduct}.
Finally, \(\sigma\) fixes the nonempty locus
\(x_1=x_2=x_3=0\); Proposition~\ref{proposition:fundamental group}
gives simple connectedness.
\end{proof}

\begin{example}
\label{example:fourfold}
The construction of Theorem~\ref{proposition:general type} has an explicit smooth model. The fixed locus
\(F=(x_1=x_2=x_3=0)\) is identified with
\(B_1\cap B_2\cap B_3\); it is smooth of codimension three.
In transverse coordinates the involution is
\((u_1,u_2,u_3)\mapsto(-u_1,-u_2,-u_3)\).
Let \(b:\widehat X=\operatorname{Bl}_F(X)\longrightarrow X\) be
the blow-up with exceptional divisor \(E\). On the chart
\((u_1,u_2/u_1,u_3/u_1)\), the lifted involution changes only
\(u_1\). The quotient variety \(Y=\widehat X/\langle\widehat\sigma\rangle\),
where \(\widehat\sigma\) denotes the lifted involution, is therefore smooth.
The quotient morphism \(q:\widehat X\longrightarrow Y\) induces a
resolution \(\rho:Y\longrightarrow Q\).
For \(D=q(E)\) with its reduced structure, \(q^*D=2E\), and the
standard blow-up and ramification formulas give
\begin{equation*}
K_{\widehat X}\sim b^*K_X+2E\sim q^*K_Y+E.
\end{equation*}
Consequently,
\begin{equation*}
K_Y\sim_{\mathbb{Q}}\rho^*K_Q+\frac{1}{2}D.
\end{equation*}
For \(n=4\), the fixed locus is a complete intersection curve of type
\((4,4,4)\) in \(\mathbb{P}^4\). Its degree is \(64\), and
adjunction gives \(\omega_F\cong\mathcal{O}_F(7)\); hence
\[
2g(F)-2=7\cdot64=448,\qquad g(F)=225.
\]
The exceptional divisor is a \(\mathbb{P}^2\)-bundle over this curve.
In particular, the 0-cycle argument does not require the entire
exceptional divisor to be rational.
\end{example}

\section{Birational type of the quotient models}
\label{section:birational-type}

For a smooth projective variety \(V\), let
\[
P_m(V)=h^0(V,\omega_V^{\otimes m})\qquad(m\geq1)
\]
be its \(m\)-th plurigenus, and let \(\kappa(V)\) be its Kodaira
dimension. Thus \(\kappa(V)=-\infty\) if all its plurigenera vanish;
otherwise, it is the largest dimension of the images of its
pluricanonical rational maps. The variety is of \emph{general type}
if \(\kappa(V)=\dim V\). For \(N=\dim V\), its canonical volume is
\[
\operatorname{vol}(K_V)=
\limsup_{m\longrightarrow\infty}\frac{N!\,P_m(V)}{m^N}.
\]
It is \emph{rationally connected} if two
general points lie on a rational curve, and it is \emph{uniruled} if
a general point lies on a rational curve. These properties, and the
Kodaira dimension, are birational invariants among smooth projective
varieties.

We use the standard vanishing of covariant holomorphic tensors on
positive dimensional rationally connected varieties and the implication
that a smooth projective uniruled variety has Kodaira dimension
\(-\infty\); see \cite{AK}. In particular, nonnegative Kodaira
dimension excludes uniruledness and rational connectedness. Vanishing
of positive-degree holomorphic forms alone does not imply rational
connectedness.

\subsection{A single multiple double cover}

\begin{proposition}\label{proposition:single-kodaira}
Use the notation and hypotheses of Theorem~\ref{theorem:single-cover},
and let \(H\subseteq G\). Define
\begin{equation}\label{equation:retained-degrees}
J_H=\{i\in[r]:\tau_i\in H\},\qquad
R_H=[r]\setminus J_H,\qquad
\delta_H=\sum_{i\in R_H}d_i.
\end{equation}
Let \(Q=X/H\), let \(\overline f:Q\longrightarrow\mathbb{P}^n\)
be the induced finite morphism, and let \(A\) be a hyperplane
divisor on \(\mathbb{P}^n\). Then \(Q\) has canonical singularities,
and
\begin{equation}\label{equation:quotient-canonical-class}
K_Q\sim_{\mathbb{Q}}(\delta_H-n-1)\overline f^{\,*}A.
\end{equation}
Every smooth projective model \(Y\) of \(Q\) satisfies
\begin{equation*}
\kappa(Y)=
\begin{cases}
-\infty,&\delta_H\leq n,\\
0,&\delta_H=n+1,\\
n,&\delta_H\geq n+2.
\end{cases}
\end{equation*}
Moreover, \(Y\) is rationally connected if and only if
\(\delta_H\leq n\); it is not uniruled in the other two cases.
In particular, it is of general type if and only if
\(\delta_H\geq n+2\).
\end{proposition}

\begin{proof}
Let \(H_0=\langle\tau_i:i\in J_H\rangle\). The covering algebra
gives
\[
X/H_0\cong X_{R_H},\qquad
Q\cong X_{R_H}/\overline H,\qquad
\overline H=H/H_0.
\]
Here \(X_{\varnothing}=\mathbb{P}^n\), as before. The group
\(\overline H\) embeds in the sign group on the retained root
coordinates and contains no coordinate involution. Indeed, a lift to
\(H\) of such an involution would differ from the corresponding
\(\tau_i\) by an element of \(H_0\subseteq H\), forcing
\(\tau_i\in H\), contrary to \(i\in R_H\).

Let \(G_{R_H}=\langle\tau_i:i\in R_H\rangle\) be the sign group
of the retained cover. At a point \(x\in X_{R_H}\), set
\(S_x=\{i\in R_H:x_i(x)=0\}\). Its stabilizer in \(G_{R_H}\) is
\(\langle\tau_i:i\in S_x\rangle\), and its stabilizer in
\(\overline H\) is the intersection with this subgroup. The simple
normal crossings hypothesis identifies the vanishing roots with part of
a regular system of local coordinates. Thus any stabilizer element
changing \(t\) roots has exactly \(t\) negative tangent eigenvalues;
all its other tangent eigenvalues are one.
Since \(\overline H\) contains no coordinate involution,
\(t\geq2\). In particular, there are no pseudoreflections, that is,
nonidentity tangent actions fixing a hyperplane pointwise.
The age, namely the sum of the normalized eigenvalue arguments in \([0,1)\),
is \(t/2\geq1\). The Reid--Shepherd-Barron--Tai criterion
\cite[Section~4.11]{Rei} implies that \(Q\) has canonical singularities.
This means that \(K_Q\) is \(\mathbb{Q}\)-Cartier and, on a
resolution \(\rho:Y\longrightarrow Q\),
\begin{equation}\label{equation:canonical-resolution}
K_Y\sim_{\mathbb{Q}}\rho^*K_Q+\sum_E a_E E,
\qquad a_E\geq0,
\end{equation}
where the sum is over the exceptional prime divisors of \(\rho\).

The same local description shows that
\(\pi:X_{R_H}\longrightarrow Q\) is \emph{quasi-\'etale}, meaning
finite and \'etale outside a subset of codimension at least two.
The ramification formula and
\eqref{equation:canonical-cover} give
\[
\pi^*K_Q\sim_{\mathbb{Q}}K_{X_{R_H}}
\sim (\delta_H-n-1)f_{R_H}^*A.
\]
Since \(f_{R_H}=\overline f\circ\pi\), proper pushforward and
division by \(\deg\pi\) prove
\eqref{equation:quotient-canonical-class}. Rational linear
equivalence is needed here: a character of the covering group can
contribute a torsion canonical class.

If \(\delta_H\leq n\), the variety \(X_{R_H}\) is smooth Fano
by Lemma~\ref{lemma:cover-geometry}, hence rationally connected
\cite{Cam,KMM}. Rational connectedness descends under a dominant
rational map, so \(Y\) is rationally connected and
\(\kappa(Y)=-\infty\).

Suppose next that \(\delta_H=n+1\). Then
\(K_Q\sim_{\mathbb{Q}}0\). For every sufficiently divisible
positive integer \(m\), canonical singularities identify
\(H^0(Y,mK_Y)\) with \(H^0(Q,mK_Q)\). The latter space is
one-dimensional. Equivalently, the effective exceptional divisor
in \eqref{equation:canonical-resolution} has no moving multiple;
see \cite[Theorem~11 and its proof]{KL}. Thus \(\kappa(Y)=0\).

Finally, if \(\delta_H\geq n+2\), the divisor \(K_Q\) is ample
by \eqref{equation:quotient-canonical-class}. Formula
\eqref{equation:canonical-resolution} expresses \(K_Y\) as a big
pullback plus an effective divisor. Hence \(K_Y\) is big and
\(\kappa(Y)=n\). In the last two cases \(Y\) is not uniruled.
All conclusions hold on every smooth projective model by birational
invariance.
\end{proof}

\Needspace{5\baselineskip}
The following proposition uses the Picard number of the canonical model to exclude birational product decompositions.

\begin{proposition}
\label{proposition:nonproduct}
Assume \(n\geq3\), and use the notation of
Proposition~\ref{proposition:single-kodaira}. The quotient
\(Q=X/H\) has Picard number one. If \(\delta_H\geq n+2\),
then \(Q\) is the canonical model of every smooth projective model
\(Y\) of \(Q\), and \(Y\) is not birational to a product of two
positive dimensional smooth projective varieties.
\end{proposition}

\begin{proof}
For a normal projective variety \(T\), write
\(N^1(T)_{\mathbb{R}}\) for the real vector space of Cartier divisors
modulo numerical equivalence. Its dimension is the Picard number.
Write \(R=R_H\). Lemma~\ref{lemma:cover-geometry} gives
\(\operatorname{Pic}(X_R)\cong\mathbb{Z}\).

The finite surjective quotient morphism \(\pi:X_R\longrightarrow Q\)
induces an injection
\[
\pi^*:N^1(Q)_{\mathbb{R}}\longrightarrow N^1(X_R)_{\mathbb{R}}.
\]
Indeed, any curve on \(Q\) is dominated by a curve on \(X_R\),
and the projection formula detects its intersection with a Cartier
divisor. Since \(Q\) has an ample divisor, its Picard number is
therefore exactly one.

If \(\delta_H\geq n+2\), Proposition~\ref{proposition:single-kodaira}
shows that \(Q\) has canonical singularities and ample canonical
\(\mathbb{Q}\)-divisor. Consequently it is the canonical model of
\(Y\). Suppose that \(Y\) were birational to \(V_1\times V_2\),
where both \(V_i\) are smooth projective and positive dimensional.
The K\"unneth formula for pluricanonical sections gives
\(P_m(V_1\times V_2)=P_m(V_1)P_m(V_2)\); since the product
is of general type, both factors are of general type. They admit
canonical models \(C_1\) and \(C_2\), by finite generation of the
canonical ring \cite{BCHM}. The product \(C_1\times C_2\) has
canonical singularities and ample canonical divisor, and is the
canonical model of \(V_1\times V_2\). Uniqueness of canonical
models gives \(Q\cong C_1\times C_2\). Pulling back an ample
Cartier divisor from each factor gives two numerically independent
classes on this product: independence follows by restricting to a
curve in each fiber of the two projections. Its Picard number is
at least two, a contradiction.
\end{proof}

\Needspace{5\baselineskip}
We next compute the fundamental group of a smooth quotient model from the elements with fixed points.

\begin{proposition}
\label{proposition:fundamental group}
Let \(n\geq2\), and retain the hypotheses of
Theorem~\ref{theorem:single-cover}. For
\(g=\tau_1^{e_1}\cdots\tau_r^{e_r}\in G\), define
\(\operatorname{wt}(g)=\sum_i e_i\), where \(e_i\in\{0,1\}\).
Define
\[
H_{\mathrm{fix}}=
\langle h\in H:\operatorname{wt}(h)\leq n\rangle.
\]
For every smooth projective model \(Y\) of \(X/H\), there is
an isomorphism
\begin{equation*}
\pi_1(Y)\cong H/H_{\mathrm{fix}}.
\end{equation*}
In particular, these fundamental groups are finite elementary abelian
\(2\)-groups.
\end{proposition}

\begin{proof}
Lemma~\ref{lemma:cover-geometry} gives \(\pi_1(X)=1\).
A covering transformation fixes a point precisely when every root
whose sign it changes vanishes there. If it changes \(t\) roots,
its fixed locus lies over the intersection of the corresponding
\(t\) branch hypersurfaces. This intersection is empty for
\(t>n\), by the simple normal crossings hypothesis. It is nonempty
for \(t\leq n\), since positive-degree hypersurfaces in
\(\mathbb{P}^n\) have nonempty intersection in that range; their
intersection has the expected dimension by transversality. The finite
cover has points over this intersection. Thus the elements of \(H\)
with fixed points are exactly those of weight at most \(n\).

Armstrong's orbit-space theorem \cite{Arm} now gives
\(\pi_1(X/H)\cong H/H_{\mathrm{fix}}\); normal closure and
subgroup generation agree because \(H\) is abelian.
Proposition~\ref{proposition:single-kodaira} shows that \(X/H\)
has canonical, hence Kawamata log terminal, singularities.
Takayama's theorem \cite[Theorem~1.1]{Tak} identifies its fundamental
group with that of a resolution. The same theorem on a common smooth
resolution gives the assertion for any smooth projective model.
\end{proof}

For a normal variety with a \(\mathbb{Q}\)-Cartier canonical divisor,
the \emph{canonical index} is the least positive integer \(a\) for
which \(aK\) is Cartier. A canonical quotient singularity is
Gorenstein precisely when its canonical divisor is Cartier.
Use the residual cover \(X_{R_H}\) and the group \(\overline H\)
from Proposition~\ref{proposition:single-kodaira}. The weight of an
element of \(\overline H\) is counted on \(R_H\), not on the
removed coordinates.

\Needspace{5\baselineskip}
The same stabilizer representations determine the terminal and Gorenstein properties of the quotient.

\begin{corollary}
\label{proposition:local-quotient-criteria}
The quotient \(Q=X/H\) is terminal if and only if \(\overline H\)
contains no element of weight two. It is Gorenstein if and only if
\(\overline H\) contains no element of odd weight at most \(n\).
More precisely, if \(\pi:X_{R_H}\longrightarrow Q\) is the quotient
map and \(Q_{\mathrm{nG}}\) denotes the non-Gorenstein locus, then
\begin{equation}\label{equation:non-gorenstein-locus}
Q_{\mathrm{nG}}=
\bigcup_{\substack{h\in\overline H\\
\operatorname{wt}(h)\text{ odd}}}
\pi\bigl(\operatorname{Fix}(h)\bigr),
\end{equation}
where \(\operatorname{Fix}(h)\) is the fixed point locus on
\(X_{R_H}\). Empty fixed point loci contribute nothing.
\end{corollary}

\begin{proof}
Apply the Reid--Shepherd-Barron--Tai criterion
\cite[Section~4.11]{Rei} and Watanabe's local invariant-ring criterion
\cite[Theorem~2]{WatII} to the stabilizer representations computed in
Proposition~\ref{proposition:single-kodaira}. There are no
pseudoreflections, the residue field is \(\mathbb{C}\), and the group
orders are invertible. An element of weight \(t\) has age \(t/2\)
and determinant \((-1)^t\), and it has a fixed point precisely when
\(t\leq n\). These observations give all the assertions.
\end{proof}

\begin{remark}\label{remark:coordinate-involutions}
The set \(J_H\) in \eqref{equation:retained-degrees} records
individual coordinate involutions belonging to \(H\). It does not
record every coordinate changed by some element of \(H\). For
example, \(H=\langle\tau_1\tau_2\tau_3\rangle\) has
\(J_H=\varnothing\), whereas \(H=\langle\tau_1\rangle\) has
\(J_H=\{1\}\). The first action has no divisorial fixed locus;
the second removes the first branch contribution from the canonical
divisor of the quotient.
\end{remark}

\begin{proposition}\label{proposition:plurigenera}
Use the notation of Proposition~\ref{proposition:single-kodaira}.
Set
\[
R=R_H,\qquad c=\delta_H-n-1,\qquad
H_0=\langle\tau_i:i\in J_H\rangle,\qquad
\overline H=H/H_0.
\]
Regard \(\overline H\) as a subgroup of the sign group on the roots
indexed by \(R\), and regard \(\chi_I\), for \(I\subseteq R\),
as a character of this sign group. For an integer \(t\), define
\[
b_n(t)=
\begin{cases}
\binom{n+t}{n},&t\geq0,\\
0,&t<0.
\end{cases}
\]
For every integer \(m\geq1\) and every smooth projective model
\(Y\) of \(Q=X/H\), one has
\begin{equation}\label{equation:plurigenera}
P_m(Y)=
\sum_{\substack{I\subseteq R\\
(\chi_I\chi_R^m)|_{\overline H}=1}}
b_n(mc-d_I).
\end{equation}
Moreover, \(2K_Q\) is Cartier and
\(2K_Q\sim2c\overline f^{\,*}A\). If \(c>0\), then
\begin{equation}\label{equation:canonical-volume}
\operatorname{vol}(K_Y)=\frac{2^r}{|H|}c^n.
\end{equation}
In particular, the varieties in
Theorem~\ref{proposition:general type} have canonical volume
four in every dimension \(n\geq4\).
\end{proposition}

\begin{proof}
By birational invariance of plurigenera, we may assume that \(Y\)
is a resolution of \(Q\) which is an isomorphism over the smooth
locus of \(Q\).
The isomorphism
\(\omega_{X_R}\cong f_R^*\mathcal{O}_{\mathbb{P}^n}(c)\)
has character twist \(\chi_R\) for the natural action on
canonical forms. To see the twist, choose local coordinates
\(z_1,\ldots,z_n\) on the base and local trivializations of the
root line bundles. The regular canonical form represented by
\[
\frac{dz_1\wedge\cdots\wedge dz_n}{\prod_{i\in R}x_i}
\]
changes by \(\chi_R(g)\) under a covering transformation \(g\).
Its regularity at the branch locus follows from the local
simple normal crossings description. The projection formula and
the covering algebra therefore give
\[
(f_R)_*\omega_{X_R}^{\otimes m}
\cong\bigoplus_{I\subseteq R}
\mathcal{O}_{\mathbb{P}^n}(mc-d_I),
\]
with character \(\chi_I\chi_R^m\) on the summand indexed by
\(I\).

By Proposition~\ref{proposition:single-kodaira}, the map
\(X_R\longrightarrow Q\) is quasi-\'etale and \(Q\) has
canonical singularities. Write
\(\omega_Q^{[m]}=\mathcal{O}_Q(mK_Q)\) for the reflexive
pluricanonical sheaf. Descent on the \'etale locus, followed by
reflexive extension across subsets of codimension at least two,
identifies
\[
H^0(Q,\omega_Q^{[m]})
\cong H^0(X_R,\omega_{X_R}^{\otimes m})^{\overline H}.
\]
Canonical singularities identify the left-hand side with
\(H^0(Y,\omega_Y^{\otimes m})\) for every \(m\geq1\).
This assertion includes the integers not divisible by the Cartier
index of \(K_Q\): take a power clearing that index. The resulting
pluricanonical section extends to a resolution, so the original
rational pluricanonical form cannot have a pole there. Conversely,
a section on the resolution restricts to the smooth locus of
\(Q\) and extends as a reflexive section. Taking invariant
summands proves \eqref{equation:plurigenera}.

For \(m=2\), the character twist of the line bundle itself is
trivial. The equivariant isomorphism
\(\omega_{X_R}^{\otimes2}\cong
f_R^*\mathcal{O}_{\mathbb{P}^n}(2c)\) descends on the
\'etale locus to
\(\omega_Q^{[2]}\cong
\overline f^{\,*}\mathcal{O}_{\mathbb{P}^n}(2c)\).
Reflexive extension gives this isomorphism on all of \(Q\),
proving the Cartier and linear-equivalence assertions.

When \(c>0\), the divisor \(K_Q\) is ample. Equality of the
pluricanonical section spaces and the Hilbert polynomial of its
Cartier multiples give
\(\operatorname{vol}(K_Y)=K_Q^n\). Since
\(\deg\overline f=2^r/|H|\), the projection formula and
\eqref{equation:quotient-canonical-class} yield
\eqref{equation:canonical-volume}. In
Theorem~\ref{proposition:general type}, one has \(c=1\),
\(r=3\), and \(|H|=2\), so the volume is four.
\end{proof}

\subsection{Character ranks and pluricanonical maps}

Whenever \(P_m(Y)>0\), write
\[
\varphi_m:Y\dashrightarrow
\mathbb{P}\bigl(H^0(Y,\omega_Y^{\otimes m})^{\vee}\bigr)
\]
for the rational map defined by the complete pluricanonical system;
\(\vee\) denotes the dual vector space. Its degree means the degree
over its image, and is used only when the map is generically finite.

In this subsection \(Y\) is a smooth projective model of a
quotient \(Q=X/H\) of general type. Use \(R,c,\overline H\) from
Proposition~\ref{proposition:plurigenera}, so \(c>0\), and write
\(\ell=\dim_{\mathbb{F}_2}(G/H)\). Here
\(\mathbb{F}_2=\mathbb{Z}/2\mathbb{Z}\). We identify subsets of \(R\)
with their characteristic vectors in \(\mathbb{F}_2^R\); their sum is
symmetric difference, denoted by \(\triangle\). For \(m\geq1\), define
\begin{equation*}
\mathcal{E}_m=\{I\subseteq R:
(\chi_I\chi_R^m)|_{\overline H}=1,\ d_I\leq mc\}.
\end{equation*}
For a nonempty \(\mathcal{E}_m\), choose \(I_0\in\mathcal{E}_m\) and define
\begin{equation*}
L_m=\operatorname{span}_{\mathbb{F}_2}
\{I\triangle I_0:I\in\mathcal{E}_m\},\qquad
\ell_m=\dim_{\mathbb{F}_2}L_m.
\end{equation*}
This subspace does not depend on the choice of \(I_0\). It lies in
the \(\ell\)-dimensional space of characters trivial on
\(\overline H\).

\Needspace{5\baselineskip}
The next theorem expresses the degree of a pluricanonical map in terms of these character differences.

\begin{theorem}
\label{theorem:character-rank}
Suppose that \(\mathcal{E}_m\) contains a subset \(J\) with
\(mc-d_J\geq1\). Then the complete \(m\)-canonical rational map
\(\varphi_m\) of \(Y\) is generically finite onto its image and
\begin{equation}\label{equation:general-map-degree}
\deg\varphi_m=2^{\ell-\ell_m}.
\end{equation}
In particular, it is birational if and only if \(\ell_m=\ell\).
The hypothesis holds for every even \(m\geq2\). Moreover,
\(\varphi_m\) is birational for every
\[
m\geq\left\lceil\frac{\delta_H+1}{c}\right\rceil.
\]
\end{theorem}

\begin{proof}
Let \(K=\mathbb{C}(\mathbb{P}^n)\) and
\(E=\mathbb{C}(X_R)\). Trivialize the root line bundles at the generic
point of the base, and write \(x_I=\prod_{i\in I}x_i\) in \(E\).
The extension \(E/K\) is Galois with group
\((\mathbb{Z}/2\mathbb{Z})^R\), and
\(\mathbb{C}(Y)=E^{\overline H}\) has degree \(2^\ell\) over \(K\).
Choose a rational generator \(\eta\) of the canonical bundle on
\(X_R\) transforming by \(\chi_R\), as in the proof of
Proposition~\ref{proposition:plurigenera}. The pullbacks of sections
from \(Y\), viewed on \(X_R\), have expressions
\(\eta^{\otimes m}\sum_{I\in\mathcal{E}_m}x_Ia_I\), where \(a_I\)
is a base polynomial of degree \(mc-d_I\). For odd \(m\), the
individual coefficients \(x_Ia_I\) need not be functions on \(Y\):
their common character under \(\overline H\) is canceled by the
character of \(\eta^{\otimes m}\).
Ratios of two sections, however, are \(\overline H\)-invariant
rational functions. On a fixed affine chart all the base polynomials
are read in the corresponding line-bundle trivializations.

Let \(K_m\) be the subfield of \(\mathbb{C}(Y)\) generated by ratios
of nonzero \(m\)-canonical sections. Ratios within the \(J\)-summand
recover \(K\), because all polynomials of a fixed positive degree
occur. Taking one nonzero section from the \(I_0\)-summand and one
from each other nonzero summand therefore gives
\[
K_m=K(x_I/x_{I_0}:I\in\mathcal{E}_m).
\]
The character of \(x_I/x_{I_0}\) is \(\chi_{I\triangle I_0}\).
Elementary Galois correspondence in the multiquadratic extension
\(E/K\) gives \([K_m:K]=2^{\ell_m}\). This proves both generic
finiteness and \eqref{equation:general-map-degree}.

For even \(m\geq2\), the empty subset lies in \(\mathcal{E}_m\)
and \(mc\geq1\). For the final assertion, \(mc\geq\delta_H+1\)
makes every subset of \(R\) satisfying the character restriction
available, with positive polynomial degree. Those subsets form an
affine coset of the full invariant-character space. Their differences
span that space, so \(\ell_m=\ell\).
\end{proof}

\begin{remark}
The positive polynomial-degree hypothesis is essential: character ranks
alone need not recover the function field of the base. For example,
the minimum volume model with \((\alpha,\beta,\gamma)=(3,3,3)\)
in Theorem~\ref{theorem:minimal-volume-maps} has dimension seven and
\(P_3(Y)=3\). The three available character blocks have constant
polynomial coefficients and their differences have full rank two,
but \(\varphi_3\) maps to \(\mathbb{P}^2\), so it is not generically
finite. Formula~\eqref{equation:general-map-degree} does not apply
in that case.
\end{remark}

\subsection{Products of multiple double covers}

\begin{proposition}\label{proposition:product-kodaira}
Use the notation and hypotheses of Theorem~\ref{theorem:products}.
Let \(\tau_{a,i}\in G_a\) be the involution changing the sign of the \(i\)-th root of \(X_a\), regarded as an element of
\(G=\prod_a G_a\) acting trivially on all factors other than
\(X_a\). For \(H\subseteq G\), define
\[
J_a=\{i\in[r_a]:\tau_{a,i}\in H\},\qquad
R_a=[r_a]\setminus J_a,\qquad
c_a=\sum_{i\in R_a}d_{a,i}-n_a-1.
\]
Let \(Q=W/H\), let \(B=\prod_a\mathbb{P}^{n_a}\), and let
\(\overline f:Q\longrightarrow B\) be the induced finite
morphism. Write \(p_a:B\longrightarrow\mathbb{P}^{n_a}\) for
the projection and \(A_a\) for a hyperplane divisor on
\(\mathbb{P}^{n_a}\). Then \(Q\) has canonical singularities and
\begin{equation}\label{equation:product-canonical-class}
K_Q\sim_{\mathbb{Q}}
\overline f^{\,*}\left(\sum_{a=1}^k c_a p_a^*A_a\right).
\end{equation}
For every smooth projective model \(Y\) of \(Q\),
\begin{equation}\label{equation:product-kodaira}
\kappa(Y)=
\begin{cases}
-\infty,&\text{some }c_a<0,\\[2pt]
\displaystyle\sum_{\{a:c_a>0\}}n_a,
&\text{all }c_a\geq0.
\end{cases}
\end{equation}
In particular, \(Y\) is of general type if and only if every
\(c_a>0\). It is uniruled if and only if some \(c_a<0\).
If all \(c_a<0\), it is rationally connected; if all
\(c_a\geq0\), it is not rationally connected.
\end{proposition}

\begin{proof}
Let \(H_0\subseteq H\) be generated by all the individual
coordinate involutions \(\tau_{a,i}\) that belong to \(H\).
Then
\[
W/H_0\cong W'=\prod_{a=1}^k X_{a,R_a},\qquad
Q\cong W'/\overline H,\qquad \overline H=H/H_0.
\]
The residual group contains no individual coordinate involution.
Every nonidentity stabilizer element therefore has at least two
negative tangent eigenvalues, counted across all factors. The
argument of Proposition~\ref{proposition:single-kodaira} proves
that \(Q\) has canonical singularities, that
\(\pi:W'\longrightarrow Q\) is quasi-\'etale, and that
\eqref{equation:product-canonical-class} holds.

If some \(c_a<0\), the corresponding partial cover is Fano and
has no nonzero pluricanonical section. The product formula for
sections of exterior tensor products gives
\[
H^0(W',mK_{W'})
=\bigotimes_{a=1}^k H^0(X_{a,R_a},mK_{X_{a,R_a}})=0
\qquad(m\geq1).
\]
Pullback by \(\pi\), together with invariance of pluricanonical
sections under a resolution of canonical singularities, implies
\(\kappa(Y)=-\infty\).

Suppose that all \(c_a\geq0\). Let
\(B_+=\prod_{\{a:c_a>0\}}\mathbb{P}^{n_a}\), and let
\(v:Q\longrightarrow B_+\) be \(\overline f\) followed by the
projection. If every \(c_a=0\), we interpret \(B_+\) as a
point. Otherwise, \eqref{equation:product-canonical-class}
expresses \(K_Q\), up to rational linear equivalence, as the
pullback by \(v\) of an ample divisor on \(B_+\). In the Stein
factorization of \(v\), the intermediate variety is finite over
\(B_+\), and the pullback of that ample divisor remains ample.
The pluricanonical image for sufficiently large and divisible
multiples therefore has dimension \(\dim B_+\). When \(B_+\)
is a point, the same conclusion follows from
\(K_Q\sim_{\mathbb{Q}}0\). Canonical singularities allow this
calculation on any smooth resolution, proving
\eqref{equation:product-kodaira}.

If some \(c_a<0\), rational curves in that Fano factor cover
\(W'\); their images under the finite map to \(Q\) are still
curves and prove uniruledness. If all \(c_a<0\), the product
\(W'\) is rationally connected, and so is its dominant rational
image \(Y\). If all \(c_a\geq0\), nonnegative Kodaira dimension
excludes uniruledness and rational connectedness.
\end{proof}

\begin{remark}
When some \(c_a<0\) and others are nonnegative, the Kodaira
dimension formula does not by itself decide rational connectedness.
Theorem~\ref{theorem:mrc} below gives a complete criterion by
recomputing the coordinate involutions after successive projections.
Independently, Theorem~\ref{theorem:products} shows that
\(s_H>0\) implies \(H^{s_H,0}(Y)\neq0\), which excludes rational
connectedness. The birational classification above is asserted for
the multiple double covers of Theorem~\ref{theorem:products}; it
does not follow from the scalar Chow action alone in
Proposition~\ref{proposition:scalar-products}.
\end{remark}

\Needspace{5\baselineskip}
The fixed point calculation also determines fundamental groups of componentwise product quotients.

\begin{corollary}
\label{corollary:product-fundamental group}
In the setting of Theorem~\ref{theorem:products}, write
\(h=(h_1,\ldots,h_k)\in H\), and let \(\operatorname{wt}(h_a)\)
be the number of root signs changed by \(h_a\) on \(X_a\). Define
\[
H_{\mathrm{fix}}^{\mathrm{prod}}
=\langle h\in H:\operatorname{wt}(h_a)\leq n_a
\text{ for every }a\rangle.
\]
Every smooth projective model \(Y\) of \(W/H\) satisfies
\[
\pi_1(Y)\cong H/H_{\mathrm{fix}}^{\mathrm{prod}}.
\]
\end{corollary}

\begin{proof}
Each factor is simply connected by Lemma~\ref{lemma:cover-geometry}.
An element fixes a point of \(W\) precisely when it fixes a point
on every factor. The fixed point calculation in
Proposition~\ref{proposition:fundamental group} gives exactly the
listed inequalities. Apply Armstrong's theorem \cite{Arm}, and then
Takayama's theorem \cite{Tak}; the quotient has canonical singularities
by Proposition~\ref{proposition:product-kodaira}.
\end{proof}

\section{Maximal rationally connected fibrations}
\label{section:mrc}

We retain the product notation of Theorem~\ref{theorem:products}.
The maximal rationally connected fibration of a smooth projective
variety is a dominant rational map whose general fibers are the
maximal rationally connected equivalence classes. Its base and the
map are unique up to birational equivalence. A fibration with
rationally connected general fiber and non-uniruled base has this
property; see the discussion preceding \cite[Corollary~1.4]{GHS}.
We now determine this fibration for every quotient under consideration.

For a nonempty subset \(S\subseteq[k]\), define
\begin{equation*}
G_S=\prod_{a\in S}G_a,\qquad
H_S=\operatorname{pr}_S(H),\qquad
Q_S=\left(\prod_{a\in S}X_a\right)/H_S,
\end{equation*}
where \(\operatorname{pr}_S:G\longrightarrow G_S\) is the
projection. For \(a\in S\), regard \(\tau_{a,i}\) as the
individual coordinate involution in \(G_S\), and define
\begin{equation*}
J_a(S)=\{i\in[r_a]:\tau_{a,i}\in H_S\},\qquad
\delta_a(S)=\sum_{i\notin J_a(S)}d_{a,i}.
\end{equation*}
Starting with \(S_0=[k]\), define recursively
\begin{equation}\label{equation:mrc-procedure}
S_{j+1}=\{a\in S_j:\delta_a(S_j)\geq n_a+1\}.
\end{equation}
This decreasing sequence stabilizes after at most \(k\) strict
decreases. Denote its stable value by \(S_\infty\). We adopt the conventions
\(G_{\varnothing}=H_{\varnothing}=\{1\}\) and
\(Q_{\varnothing}=\operatorname{Spec}\mathbb{C}\), and keep the
empty set unchanged by the recursion.

\begin{theorem}\label{theorem:mrc}
Let \(Y\) be a smooth projective model of \(W/H\), and let
\(Y_\infty\) be a smooth projective model of \(Q_{S_\infty}\),
with \(Y_\infty=\operatorname{Spec}\mathbb{C}\) when
\(S_\infty=\varnothing\). The natural projection induces the
maximal rationally connected fibration
\[
Y\dashrightarrow Y_\infty.
\]
In particular,
\begin{equation}\label{equation:mrc-criterion}
Y\text{ is rationally connected}
\quad\Longleftrightarrow\quad S_\infty=\varnothing.
\end{equation}
The dimension and Kodaira dimension of its base are
\begin{equation}\label{equation:mrc-base-invariants}
\dim Y_\infty=\sum_{a\in S_\infty}n_a,
\qquad
\kappa(Y_\infty)=
\sum_{\substack{a\in S_\infty\\
\delta_a(S_\infty)\geq n_a+2}}n_a.
\end{equation}
The Kodaira dimension of a point is taken to be zero.
Equivalently, \(Y\) is rationally connected if and only if every
nonempty subset \(S\subseteq[k]\) contains an index \(a\) with
\(\delta_a(S)\leq n_a\).
\end{theorem}

\begin{proof}
Consider one step of \eqref{equation:mrc-procedure}, with current
set \(S\). Define
\[
D=\{a\in S:\delta_a(S)\leq n_a\},\qquad T=S\setminus D,
\]
and suppose that \(D\neq\varnothing\). Let
\[
H_0(S)=\langle\tau_{a,i}:a\in S,\ i\in J_a(S)\rangle,
\qquad \overline H_S=H_S/H_0(S).
\]
First quotient by \(H_0(S)\), as in
Proposition~\ref{proposition:product-kodaira}. Writing
\(R_a(S)=[r_a]\setminus J_a(S)\), one obtains
\[
Q_S\cong(F\times V)/\overline H_S,
\qquad
F=\prod_{a\in D}X_{a,R_a(S)},\qquad
V=\prod_{a\in T}X_{a,R_a(S)}.
\]
Each factor of \(F\) is smooth Fano by
Lemma~\ref{lemma:cover-geometry}, because its retained degree sum
is at most \(n_a\). Thus \(F\) is rationally connected.

Let \(L\) be the image of \(\overline H_S\) in the group of
automorphisms of \(V\), and let
\(K=\ker(\overline H_S\longrightarrow L)\). Quotienting the
retained factors shows that \(V/L\cong Q_T\). Choose a dense
open subset of \(V\) on which \(L\) acts freely. Over its image,
and after choosing a lift of a geometric point, the fiber of
\(Q_S\longrightarrow Q_T\) is isomorphic to \(F/K\).
Indeed, the lifts of that point form an \(L\)-orbit, and the
subgroup preserving a chosen lift is precisely \(K\).
This description also applies when \(T=\varnothing\), with
\(V\) a point. The geometric general fiber is therefore integral and has a smooth
projective model dominated rationally by \(F\). Rational connectedness
descends under this dominant rational map, so that model is rationally
connected.

Choose compatible smooth projective models of this finite sequence of
quotients and resolve the projection maps. After shrinking the bases,
generic smoothness and the preceding fiber description identify their
geometric general fibers with smooth projective models of the indicated
quotients. Each such fiber is rationally connected. The composite has
rationally connected geometric general fiber as well: apply
\cite[Corollary~1.3]{GHS} successively to the fiber tower over an
algebraic closure of the function field of the final base. If
\(S_\infty=\varnothing\), this proves rational connectedness of \(Y\).

If \(S_\infty\neq\varnothing\), stabilization means that
\(\delta_a(S_\infty)\geq n_a+1\) for every remaining index.
Proposition~\ref{proposition:product-kodaira} proves
\eqref{equation:mrc-base-invariants} and shows that
\(Y_\infty\) is not uniruled. The standard characterization of the
maximal rationally connected fibration by rationally connected general
fiber and non-uniruled base, recalled before the theorem, now applies.
It gives maximality and \eqref{equation:mrc-criterion}.

For the final criterion, observe that \(S\subseteq T\) and
\(a\in S\) imply
\[
J_a(S)\supseteq J_a(T),\qquad
\delta_a(S)\leq\delta_a(T).
\]
Call a subset \(S\) stable if
\(\delta_a(S)\geq n_a+1\) for every \(a\in S\).
The displayed monotonicity shows inductively that every stable
subset is contained in every \(S_j\). Since \(S_\infty\)
is itself stable, it is the largest stable subset. It is empty
exactly when no nonempty stable subset exists, as asserted.
\end{proof}

\begin{corollary}\label{corollary:mrc-support}
With the notation of Theorem~\ref{theorem:mrc}, one has
\[
s_0(Y)=s_0(Y_\infty),\qquad
h^{p,0}(Y)=h^{p,0}(Y_\infty)\quad(p\geq0).
\]
\end{corollary}

\begin{proof}
Consider one step \(S\longmapsto T\), and use the admissible
tuples from \eqref{equation:admissible-factor} and
\eqref{equation:product-bound} for the factors indexed by \(S\).
If \(\boldsymbol{I}=(I_a)_{a\in S}\in\mathcal{A}(H_S)\),
invariance under \(\tau_{a,i}\in H_S\) gives
\(i\notin I_a\) for every \(i\in J_a(S)\). Hence
\(d_{a,I_a}\leq\delta_a(S)\). For \(a\in S\setminus T\),
the right-hand side is at most \(n_a\), so admissibility forces
\(I_a=\varnothing\). Restriction to \(T\) therefore gives a
bijection
\(\mathcal{A}(H_S)\longrightarrow\mathcal{A}(H_T)\), with
inverse extension by empty subsets. It preserves both weights and
multiplicities, because \(H_T=\operatorname{pr}_T(H_S)\).
Theorem~\ref{theorem:products} proves both assertions at each step.
If the terminal set is empty, its unique empty tuple has weight
zero and multiplicity one, in agreement with the assertions for
a point.
\end{proof}

\Needspace{5\baselineskip}
The next two examples illustrate the effect of recomputing the projected subgroup in mixed cases.

\begin{example}
\label{example:mixed-rc}
Let \(F\longrightarrow\mathbb{P}^2\) be the double cover
branched along a smooth conic, with covering involution \(\alpha\).
It is a smooth quadric surface, isomorphic to
\(\mathbb{P}^1\times\mathbb{P}^1\). Let
\(T\longrightarrow\mathbb{P}^2\) be a double cover branched
along a smooth sextic, with covering involution \(\beta\).
The surface \(T\) is a K3 surface. Take
\(H=\langle(\alpha,\beta)\rangle\) on \(F\times T\).
There is no individual coordinate involution in \(H\), so the
initial canonical coefficients in
Proposition~\ref{proposition:product-kodaira} are \((-2,0)\).
The first step of \eqref{equation:mrc-procedure} removes \(F\).
The projected subgroup on \(T\) is then \(\langle\beta\rangle\),
which removes its only root and gives \(T/\langle\beta\rangle
\cong\mathbb{P}^2\). The second step therefore leaves a point.
Every smooth projective quotient model is rationally connected,
has \(\kappa=-\infty\), and has trivial degree zero Chow group. In particular, it is not of general type.
\end{example}

\Needspace{5\baselineskip}
In contrast, the next fourfold is uniruled but not rationally connected.

\begin{example}
\label{example:mixed-enriques}
Let \(F\) and \(\alpha\) be as in
Example~\ref{example:mixed-rc}. Let \(S\) be the K3 surface in
Example~\ref{example:hu}, formed from three conic branches, and
let \(\epsilon=\tau_1\tau_2\tau_3\). This involution is free,
and \(E=S/\langle\epsilon\rangle\) is an Enriques surface.
The diagonal action of \(\langle(\alpha,\epsilon)\rangle\)
on \(F\times S\) is free, so
\[
Y=(F\times S)/\langle(\alpha,\epsilon)\rangle
\]
is a smooth projective fourfold. Its initial canonical coefficients
are again \((-2,0)\). The first step removes \(F\), but the
projected subgroup \(\langle\epsilon\rangle\) contains no
individual coordinate involution on \(S\). The procedure stops
with base \(E\). The morphism
\(Y\longrightarrow E\) has every geometric fiber isomorphic to
\(F\) and is its maximal rationally connected fibration.
Thus \(Y\) is uniruled but not rationally connected, and
\(\kappa(Y)=-\infty\), so it is not of general type.
Nevertheless, \(s_0(E)=0\) by Theorem~\ref{theorem:single-cover},
and Corollary~\ref{corollary:mrc-support} gives
\(s_0(Y)=0\). Hence the degree map on
\(\operatorname{CH}_0(Y)\) is an isomorphism, and all its
positive-degree holomorphic forms vanish. These two examples show
why the projected subgroup must be recomputed at every step.
\end{example}

\section{Sharp bounds and point-support models}
\label{section:branch-bound}

We return to a single multiple double cover, with the notation of
Theorem~\ref{theorem:single-cover} and
Proposition~\ref{proposition:single-kodaira}. Identify the sign group
\(G\) with \(\mathbb{F}_2^r\) using its coordinate involutions,
where \(\mathbb{F}_2=\mathbb{Z}/2\mathbb{Z}\). If
\(n=\sum_{e\geq0}\nu_e2^e\), with \(\nu_e\in\{0,1\}\),
write \(s_2(n)=\sum_e\nu_e\) for the sum of its binary digits.

\begin{theorem}\label{theorem:branch-bound}
Suppose that the degree map on \(\operatorname{CH}_0(Y)\) is
an isomorphism for a smooth projective model \(Y\) of \(X/H\).
Let \(\ell=\dim_{\mathbb{F}_2}(G/H)\). Then
\begin{equation}\label{equation:sharp-branch-bound}
\delta_H\leq
\sum_{j=0}^{\ell-1}\left\lfloor\frac{n}{2^j}\right\rfloor
\leq 2n-s_2(n).
\end{equation}
For \(\ell=0\), the first sum is zero. The final bound is attained
for every integer \(n\geq2\) by a cover with all \(d_i=1\)
and a subgroup containing no individual coordinate involution.
\end{theorem}

\begin{proof}
Let
\[
C=H^{\perp}=\{v\in\mathbb{F}_2^r:
v\cdot h=0\text{ for every }h\in H\},
\]
where \(v\cdot h=\sum_i v_i h_i\) in \(\mathbb{F}_2\).
This is a binary linear code, meaning a vector subspace of
\(\mathbb{F}_2^r\), of dimension \(\ell\). Its words are the
characteristic vectors of the subsets in \(\mathcal{I}(H)\).
A coordinate vanishes on every word precisely when the corresponding
coordinate vector belongs to \(C^{\perp}=H\). Thus deleting all
zero coordinates deletes exactly \(J_H\) and leaves a code with
full support on \(R_H\). Here full support means that every
remaining coordinate is nonzero on some word.

Repeat the coordinate indexed by \(i\in R_H\) exactly \(d_i\)
times. This gives a binary linear code \(C'\) of length
\(\delta_H\), dimension \(\ell\), and full support. The
Hamming weight of a word, namely the number of its nonzero
coordinates, is \(d_I\) for the corresponding invariant subset
\(I\). By Theorem~\ref{theorem:single-cover}, every such weight
is at most \(n\).

The binary antiGriesmer bound states that a full-support linear code
of dimension \(\ell\) and maximum weight \(w\) has length at most
\(\sum_{j=0}^{\ell-1}\lfloor w/2^j\rfloor\);
see \cite[Theorem~1.1]{ZCLL}. Its hypothesis is full support;
repeated coordinates, as used here, are allowed. Applying this bound
to \(C'\),
with \(w\leq n\), proves the first inequality in
\eqref{equation:sharp-branch-bound}. The identity
\[
\sum_{j\geq0}\left\lfloor\frac{n}{2^j}\right\rfloor
=\sum_{e\geq0}\nu_e(2^{e+1}-1)=2n-s_2(n)
\]
proves the second.

We prove sharpness by the standard binary simplex-code construction.
For each exponent \(e\) with \(\nu_e=1\), take the code
\[
C_e=\{(a\cdot v)_{v\in\mathbb{F}_2^{e+1}\setminus\{0\}}:
a\in\mathbb{F}_2^{e+1}\}.
\]
It has length \(2^{e+1}-1\), dimension \(e+1\), and full
support. Every nonzero word has weight \(2^e\), because a
nonzero linear functional takes the value one on exactly half
of \(\mathbb{F}_2^{e+1}\). The direct sum
\(C=\bigoplus_{\{e:\nu_e=1\}}C_e\) therefore has full support,
length
\(r=2n-s_2(n)\), and maximum weight \(n\).

Choose \(r\) general quadric hypersurfaces in \(\mathbb{P}^n\).
They satisfy \eqref{equation:branch-hypotheses}. Form the associated
multiple double cover, with \(d_i=1\), and choose \(H=C^{\perp}\)
in its sign group. Full support of \(C\) implies
\(J_H=\varnothing\), so \(\delta_H=r=2n-s_2(n)\).
The maximum-weight condition and
Theorem~\ref{theorem:single-cover} give
\(\operatorname{CH}_0(Y)\cong\mathbb{Z}\) by the degree map.
This proves sharpness. The combinatorial bound and simplex codes
are standard coding-theoretic ingredients; their application here
concerns the branch degree and 0-cycles of these quotient models.
\end{proof}

\begin{corollary}\label{corollary:dimension-threshold}
Among the quotients in Theorem~\ref{theorem:single-cover}, a smooth
projective model of dimension \(n\) that is of general type and has
trivial degree zero Chow group exists if and only if \(n\geq4\).
\end{corollary}

\begin{proof}
For \(n=2\), Theorem~\ref{theorem:branch-bound} gives
\(\delta_H\leq3\); for \(n=3\), it gives
\(\delta_H\leq4\). General type would require
\(\delta_H\geq n+2\) by
Proposition~\ref{proposition:single-kodaira}, which is impossible
in either case. Theorem~\ref{proposition:general type} supplies
examples for every \(n\geq4\). The sharpness construction in
Theorem~\ref{theorem:branch-bound} supplies another family:
for \(n\geq4\),
\[
\delta_H-n-1=n-s_2(n)-1\geq2.
\]
Indeed, the binary expansion then has a nonzero digit of exponent
at least two. For \(n=2,3\), that construction instead has
\(\delta_H=n+1\) and Kodaira dimension zero.
\end{proof}

The dimension restriction in this corollary uses both the independent
double-cover construction and the smooth simple normal crossings
branch hypotheses. It is not a restriction on arbitrary varieties
with \(\operatorname{CH}_0\cong\mathbb{Z}\); in particular,
the threefold of general type in \cite[Example~19]{Pet} is obtained
by a different construction.

\subsection{Fundamental groups when the degree map is an isomorphism}

\begin{theorem}\label{theorem:fundamental-classification}
Let \(Y\) be a smooth projective model of \(X/H\) under the hypotheses of
Theorem~\ref{theorem:single-cover}, and suppose that
\(\deg:\operatorname{CH}_0(Y)\longrightarrow\mathbb{Z}\)
is an isomorphism. Then \(Y\) is simply connected, except in the
following case. After quotienting by the individual coordinate
involutions in \(H\), the retained cover has exactly \(n+1\)
roots, each of degree parameter one, the integer \(n\) is even,
and the residual subgroup is generated by the simultaneous change
of all retained root signs. In that case
\[
\pi_1(Y)\cong\mathbb{Z}/2\mathbb{Z},\qquad \kappa(Y)=0.
\]
In particular, every general type quotient model in this covering
class with \(A_0(Y)=0\) is simply connected.
\end{theorem}

\begin{proof}
Remove the coordinate involutions as in
Proposition~\ref{proposition:single-kodaira}. Relabel the retained
roots by \(1,\ldots,s\), and write \(L\subseteq\mathbb{F}_2^s\)
for the residual subgroup. No coordinate vector belongs to \(L\).
The code \(C=L^{\perp}\) consequently has full support. Write
\(\ell=\dim_{\mathbb{F}_2}C\). Its words have weighted degree at
most \(n\), by Theorem~\ref{theorem:single-cover}. If \(\ell=0\),
full support forces \(s=0\), and the quotient is projective space.
We may therefore assume \(\ell>0\).

Choose a set of \(\ell\) independent columns of a generator
matrix of \(C\), and permute the coordinates so that the matrix
has systematic form
\[
\bigl(I_{\ell}\mid M\bigr).
\]
Here \(I_{\ell}\) is the identity matrix. Write
\(t=s-\ell\) and let \(v_1,\ldots,v_t\) be the columns of
\(M\); each \(v_j\) is nonzero. The sum of all rows gives a
word whose first \(\ell\) coordinates are one. Positivity of
the degree parameters therefore implies \(\ell\leq n\).
The subgroup \(L\) has the basis
\[
h_j=(v_j,e_j)\in\mathbb{F}_2^{\ell}\oplus\mathbb{F}_2^t,
\qquad 1\leq j\leq t,
\]
where \(e_j\) is the \(j\)-th coordinate vector. If \(\ell<n\),
each basis vector has weight at most \(\ell+1\leq n\).
Proposition~\ref{proposition:fundamental group} then gives
\(\pi_1(Y)=1\).

It remains to consider \(\ell=n\). The sum of all rows already
has weighted degree at least \(n\) in its first \(n\)
coordinates. Since its total weighted degree is at most \(n\),
those first \(n\) degree parameters equal one, and every column
\(v_j\) has even weight. For each \(i\in\{1,\ldots,n\}\),
use the word obtained by summing all rows except row \(i\).
Its first \(n\) coordinates contribute \(n-1\). On column
\(v_j\), its value is \((v_j)_i\), because \(v_j\) has
even weight. Hence
\begin{equation*}
\sum_{\{j:(v_j)_i=1\}}d_{n+j}\leq1
\qquad(1\leq i\leq n).
\end{equation*}
Every \(v_j\) is nonzero. Thus all the remaining degree
parameters also equal one, and the supports of the columns
\(v_j\) are pairwise disjoint.

The basis vector \(h_j\) has weight
\(1+\operatorname{wt}(v_j)\). Unless some \(v_j\) has weight
\(n\), all these weights are at most \(n\), and the
fundamental group is again trivial. If one column has weight
\(n\), its even weight forces \(n\) to be even. Its support
uses every row, so disjointness and nonvanishing exclude every
other column. Therefore \(s=n+1\), all \(d_i=1\), and
\(L=\langle(1,\ldots,1)\rangle\). Its nonidentity element
has no fixed point, and
Proposition~\ref{proposition:fundamental group} gives
\(\pi_1(Y)\cong\mathbb{Z}/2\mathbb{Z}\).
The retained degree sum is \(n+1\), so
Proposition~\ref{proposition:single-kodaira} gives \(\kappa(Y)=0\).
Conversely, this exceptional configuration is the free quotient in
Example~\ref{example:hu}; it has the asserted Chow group and
fundamental group. This also proves the final assertion.
\end{proof}

\subsection{Gorensteinness and the canonical index}

\begin{theorem}
\label{theorem:gorenstein-obstruction}
Suppose that \(A_0(Y)=0\) for a smooth projective model of \(Q=X/H\) as in
Theorem~\ref{theorem:single-cover}, and that \(\kappa(Y)\geq0\).
Then \(Q\) is Gorenstein if and only if the residual cover and
subgroup are the exceptional even-dimensional free quotient in
Theorem~\ref{theorem:fundamental-classification}. In every other
case, \(Q\) has canonical index exactly two and has a nonempty
non-Gorenstein locus. In particular, every general type point-support
canonical model in this covering class is non-Gorenstein of index two.
\end{theorem}

\begin{proof}
Write \(R=R_H\) and \(L=\overline H\). The assumption on Kodaira
dimension implies \(\delta_H\geq n+1\). The Chow criterion in
Theorem~\ref{theorem:single-cover} consequently excludes
\(\chi_R|_L=1\), because the full retained subset has degree greater
than \(n\). Thus the total sign character is nontrivial on \(L\).

If \(Q\) is Gorenstein, Corollary~\ref{proposition:local-quotient-criteria}
shows that every element of \(L\) with a fixed point has even weight.
The subgroup they generate lies in the proper subgroup
\(\ker(\chi_R|_L)\). Proposition~\ref{proposition:fundamental group}
therefore gives \(\pi_1(Y)\ne1\). Theorem~\ref{theorem:fundamental-classification}
identifies this with the specified exceptional case. Conversely,
that residual action is free, so its quotient is smooth and hence
Gorenstein. In all other cases \(K_Q\) is not Cartier, whereas
\(2K_Q\) is Cartier by Proposition~\ref{proposition:plurigenera}.
Its index is therefore two. The exceptional case has Kodaira
dimension zero, proving the general type assertion.
\end{proof}

\begin{remark}
Canonical index is not the order of the canonical class in the Picard
group. In the exceptional smooth quotient, the canonical index is one,
although its canonical line bundle has order two. The index-two
assertion above concerns the singular canonical model, not the smooth
resolution.
\end{remark}

\subsection{Minimum canonical volume and the equality case}

\begin{theorem}\label{theorem:minimum volume}
Let \(Y\) be a smooth projective model of \(X/H\) under the
hypotheses of Theorem~\ref{theorem:single-cover}. Suppose that \(Y\) is of
general type and that its degree map on 0-cycles is an
isomorphism. Then
\begin{equation}\label{equation:minimum-canonical-volume}
\operatorname{vol}(K_Y)\geq4.
\end{equation}
Equality holds if and only if
\begin{equation}\label{equation:minimum volume-equality}
\dim_{\mathbb{F}_2}(G/H)=2,
\qquad \delta_H=n+2.
\end{equation}

The equality case has the following equivalent character description.
On the retained root set \(R_H\), choose a basis of the two-dimensional
space of \(H\)-invariant characters. Partition \(R_H\) into
\(R_1,R_2,R_3\) according as the corresponding column of this
basis is \((1,0)\), \((0,1)\), or \((1,1)\), and define
\(\lambda_j=\sum_{i\in R_j}d_i\). Then
\begin{equation*}
\lambda_1+\lambda_2+\lambda_3=n+2,
\qquad \lambda_j\geq2\quad(j=1,2,3).
\end{equation*}
Conversely, every partition with these properties and a
two-dimensional invariant-character space gives a general type
quotient with \(\operatorname{CH}_0(Y)\cong\mathbb{Z}\) and
canonical volume four. In particular, every triple of integers
\(\lambda_1,\lambda_2,\lambda_3\geq2\) with sum \(n+2\)
is realized by three smooth branch hypersurfaces of degrees
\(2\lambda_1,2\lambda_2,2\lambda_3\) and the diagonal involution.
\end{theorem}

\begin{proof}
Write \(\ell=\dim_{\mathbb{F}_2}(G/H)\) and
\(c=\delta_H-n-1\). General type means \(c\geq1\).
If \(\ell=0\), there are no retained roots, so \(c<0\).
If \(\ell=1\), the unique nonzero invariant character is nonzero
on every retained coordinate. Its weighted degree is therefore
\(\delta_H\). The criterion of Theorem~\ref{theorem:single-cover}
would imply \(\delta_H\leq n\), again contradicting \(c\geq1\).
Consequently \(\ell\geq2\). The volume formula
\eqref{equation:canonical-volume} gives
\[
\operatorname{vol}(K_Y)=2^{\ell}c^n\geq4.
\]
Since \(\ell\) and \(c\) are integers, equality is equivalent
to \(\ell=2\) and \(c=1\), which is
\eqref{equation:minimum volume-equality}.

There are no zero columns on \(R_H\), by its definition. When
\(\ell=2\), its nonzero column types are exactly the three listed
in the statement. The three nonzero invariant characters have
weighted degrees
\[
\lambda_1+\lambda_3,\qquad
\lambda_2+\lambda_3,\qquad
\lambda_1+\lambda_2.
\]
Since \(\lambda_1+\lambda_2+\lambda_3=\delta_H=n+2\), the condition that each
of these three degrees is at most \(n\) is exactly
\(\lambda_j\geq2\) for every \(j\). This proves both the
necessity and the converse.

For the last assertion, choose three general hypersurfaces of the
indicated degrees in \(\mathbb{P}^n\), and use their independent
double covers. The diagonal subgroup is generated by
\((1,1,1)\), and its invariant code consists of the even subsets
of three roots. The three nonzero weights are
\(n+2-\lambda_j\leq n\). The Chow criterion, the canonical
class formula, and the volume formula prove the assertion.
These equality models are simply connected by
Theorem~\ref{theorem:fundamental-classification}, and they are not
birational products by Proposition~\ref{proposition:nonproduct}.
\end{proof}

\begin{remark}
The minimum in \eqref{equation:minimum-canonical-volume} is a minimum
inside the covering class of this paper, not among all varieties
with \(A_0=0\). For example, Vial's fourfold \cite[Section~2.3]{Via}
is an \'etale quotient of degree seven of the degree-seven Fermat
hypersurface in \(\mathbb{P}^5\). Adjunction and the degree
formula give its canonical volume as \(7/7=1\).
\end{remark}

\Needspace{5\baselineskip}
We next determine the singularities of the canonical models in the equality case.

\begin{corollary}
\label{corollary:minimal-volume-singularities}
In the equality case of Theorem~\ref{theorem:minimum volume}, the
canonical model \(Q\) is terminal if and only if each of the three
character blocks \(R_1,R_2,R_3\) consists of one root. In every equality
case, the non-Gorenstein locus has codimension three, and its general
transverse singularity is \(\frac{1}{2}(1,1,1)\), the quotient of
\(\mathbb{A}^3\) by scalar multiplication by \(-1\).
\end{corollary}

\begin{proof}
The residual subgroup is the kernel of a two-row matrix whose columns
are the three distinct nonzero vectors of \(\mathbb{F}_2^2\). It
contains a weight-two vector precisely when two columns are equal,
that is, when one block has at least two roots. The terminal criterion
is Corollary~\ref{proposition:local-quotient-criteria}.

For a kernel vector, the parities of the numbers of chosen coordinates
in the three blocks are equal. Such a vector has odd weight precisely
when all three parities are one. Its support then contains a triple
with one coordinate from each block, itself a kernel vector. Its fixed
locus is contained in the fixed locus of that triple. Formula
\eqref{equation:non-gorenstein-locus} therefore expresses the
non-Gorenstein locus as the union of the images of these triple fixed
loci. Each has codimension three and is nonempty. At a general point
on one such locus no other root vanishes, so the stabilizer has order
two and acts by minus the identity in the three normal directions.
\end{proof}

\Needspace{5\baselineskip}
The three-block description also determines the complete pluricanonical systems of every minimum volume model.

\begin{theorem}
\label{theorem:minimal-volume-maps}
Suppose that \(Y\) is in the equality case of
Theorem~\ref{theorem:minimum volume}. Relabel the three character
blocks so that their degree sums satisfy
\[
2\leq\alpha\leq\beta\leq\gamma,
\qquad \alpha+\beta+\gamma=n+2.
\]
Use \(b_n\) from Proposition~\ref{proposition:plurigenera}. Then
\begin{equation}\label{equation:all-minimal-plurigenera}
P_m(Y)=
\begin{cases}
b_n(m)+b_n(m-\alpha-\beta)+b_n(m-\alpha-\gamma)
+b_n(m-\beta-\gamma),&m\text{ even},\\
b_n(m-\alpha)+b_n(m-\beta)+b_n(m-\gamma)
+b_n(m-n-2),&m\text{ odd}.
\end{cases}
\end{equation}
For even \(m\geq2\), the complete pluricanonical map is generically
finite and
\begin{equation}\label{equation:minimal-even-degree}
\deg\varphi_m=
\begin{cases}
4,&m<\alpha+\beta,\\
2,&\alpha+\beta\leq m<\alpha+\gamma,\\
1,&m\geq\alpha+\gamma.
\end{cases}
\end{equation}
For odd \(m\geq\alpha+1\), it is generically finite and
\begin{equation}\label{equation:minimal-odd-degree}
\deg\varphi_m=
\begin{cases}
4,&m<\beta,\\
2,&\beta\leq m<\gamma,\\
1,&m\geq\gamma.
\end{cases}
\end{equation}
For odd \(m<\alpha\) there are no sections. If \(\alpha\) is odd,
the \(\alpha\)-canonical map has image of dimension at most two and
is not generically finite.

Let \(t=\max\{\gamma,\alpha+1\}\). The first birational
pluricanonical map occurs at
\begin{equation}\label{equation:minimal-first-birational}
2\left\lceil\frac{t-1}{2}\right\rceil+1,
\end{equation}
and the least integer \(m_0\) such that every \(\varphi_m\),
\(m\geq m_0\), is birational is
\begin{equation}\label{equation:minimal-eventual-birational}
m_0=2\left\lceil\frac{\alpha+\gamma}{2}\right\rceil-1.
\end{equation}
\end{theorem}

\begin{proof}
The four invariant subsets of the retained root set are the empty
subset and the three unions of two character blocks. Since \(c=1\),
these give the even formula in \eqref{equation:all-minimal-plurigenera}.
For odd \(m\), the character condition translates this set by the
full retained set, giving the three single blocks and the full set.
Their degree sums are \(\alpha,\beta,\gamma,n+2\), proving the odd
formula.

For even \(m\geq2\), use the empty subset as reference in
Theorem~\ref{theorem:character-rank}. The first nontrivial character
appears at \(\alpha+\beta\); two independent characters appear
at \(\alpha+\gamma\). This proves
\eqref{equation:minimal-even-degree}. For odd \(m\geq\alpha+1\),
use the block of degree \(\alpha\) as reference. Its positive
polynomial degree recovers the base field. The other blocks appear
at \(\beta\) and \(\gamma\), contributing successively one and
two independent character differences. The full subset appears later
and does not change the span. This proves
\eqref{equation:minimal-odd-degree}.

For odd \(m<\alpha\), every summand is zero. When \(m=\alpha\)
is odd, only blocks of degree \(\alpha\) contribute, each a
one-dimensional space. There are at most three such blocks, so the
image has dimension at most two, less than \(n\geq4\).

The first odd birational degree is the least odd integer at least
\(t\), giving \eqref{equation:minimal-first-birational}. The first
even one is the least even integer at least \(\alpha+\gamma\).
Because \(\alpha\geq2\), the former is smaller than the latter.
Within each parity, the displayed degree formulas are nonincreasing.
The largest nonbirational even degree is two less than the first
birational even degree, and every larger odd degree is already
birational. This proves \eqref{equation:minimal-eventual-birational}.
\end{proof}

\begin{corollary}
\label{theorem:pluricanonical maps}
For the varieties of Theorem~\ref{proposition:general type}, take
\((\alpha,\beta,\gamma)=(2,2,n-2)\) in
Theorem~\ref{theorem:minimal-volume-maps}. All maps \(\varphi_m\)
with \(m\geq2\) are generically finite. Their first birational degree
and their eventual birationality threshold are, respectively,
\[
2\lfloor n/2\rfloor-1
\qquad\text{and}\qquad
2\lceil n/2\rceil-1.
\]
\end{corollary}

\begin{proof}
Substitute the indicated triple into
\eqref{equation:minimal-even-degree},
\eqref{equation:minimal-odd-degree},
\eqref{equation:minimal-first-birational}, and
\eqref{equation:minimal-eventual-birational}. The positive-degree
hypothesis holds for every even \(m\geq2\) and every odd \(m\geq3\).
\end{proof}

For example, when \(n=5\), the map \(\varphi_3\) is birational,
\(\varphi_4\) has degree two, and all \(\varphi_m\) with \(m\geq5\)
are birational. The first birational degree must therefore be
distinguished from the eventual threshold.

\begin{example}
In dimension ten, the triples \((2,2,8)\) and \((4,4,4)\) both give
simply connected general type models with \(A_0=0\), canonical volume
four, and no birational product decomposition. For \((2,2,8)\), the
first and eventual birational degrees are both nine. For \((4,4,4)\),
they are five and seven, respectively. Thus neither the canonical
volume nor the vanishing of \(A_0\) determines these thresholds.
\end{example}

\Needspace{5\baselineskip}
We finish with a fourfold defined by seven quadric branch hypersurfaces.

\begin{example}
\label{example:seven-quadrics}
For \(n=4\), the sharpness construction uses the binary simplex
code of length seven and dimension three; all its nonzero words
have weight four. Take seven general quadric branches in
\(\mathbb{P}^4\) and let \(H\) be the orthogonal complement of
this code in \(\mathbb{F}_2^7\). Then \(|H|=16\),
\(J_H=\varnothing\), and \(\delta_H=7\). The quotient
\(Q=X/H\) is finite of degree eight over \(\mathbb{P}^4\),
has canonical singularities, and satisfies
\[
K_Q\sim_{\mathbb{Q}}2\overline f^{\,*}A.
\]
Every smooth projective model \(Y\) is of general type, is not
uniruled or rationally connected, and has
\(\operatorname{CH}_0(Y)\cong\mathbb{Z}\). It is not birational
to a product, by Proposition~\ref{proposition:nonproduct}. By
Proposition~\ref{proposition:plurigenera},
\(\operatorname{vol}(K_Y)=8\cdot2^4=128\). This is the
canonical volume, not an assertion about the top self-intersection
of the canonical divisor on an arbitrary resolution.
\end{example}

\section*{AI declaration}
AI-assisted tools (ChatGPT 6) were used for literature searches, proof checking, and
preparation of the manuscript.
Responsibility for the mathematical content and references
rests with the author.

\end{document}